\documentclass[12pt]{article}
\usepackage{amsmath}
\usepackage{amssymb}

\usepackage{graphicx}

\def\qed{\hfill$\Box$\par\medskip\par\relax}

\def\1#1{{\bf 1}{\{#1\}}}

\newcommand{\eps}{\varepsilon}
\newcommand{\Z}{{\mathbb Z}}

\newcommand{\Sph}{{\mathbb S}}

\newcommand{\HH}{{\mathcal H}}
\newcommand{\LL}{{\mathcal L}}
\newcommand{\R}{{\mathbb R}}
\newcommand{\RR}{{\mathcal R}}
\newcommand{\tRR}{{\tilde{\mathcal R}}}
\newcommand{\CC}{{\mathfrak C}}
\newcommand{\cCC}{{\mathcal C}}
\newcommand{\OO}{{\mathcal O}}

\newcommand{\DD}{{\mathcal D}}
\newcommand{\BB}{{\mathcal B}}
\newcommand{\cDD}{\bar{\mathcal D}}
\newcommand{\FF}{{\mathfrak F}}

\let\phi=\varphi

\newcommand{\PP}{{\mathbf P}}
\newcommand{\E}{{\mathbf E}}

\newcommand{\ttau}{{\tilde\tau}}
\newcommand{\n}{{\mathbf n}}

\newcommand{\m}{{\mathfrak m}}
\newcommand{\bm}{{\bar{\mathfrak m}}}

\newcommand{\bga}{{\bar\gamma}}
\newcommand{\tX}{{\tilde X}}
\newcommand{\tV}{{\tilde V}}
\newcommand{\tK}{{\tilde K}}
\newcommand{\txi}{{\tilde\xi}}
\newcommand{\htau}{{\hat\tau}}

\newcommand{\fr}{\partial\DD}

\newcommand{\diam}{{\mathop{\rm diam}}}
\newcommand{\8}{{\infty}}

\newcommand{\eqlaw}{\stackrel{\text{\tiny law}}{=}}

\allowdisplaybreaks

\newtheorem{theo}{Theorem}[section]
\newtheorem{lmm}{Lemma}[section]
\newtheorem{df}{Definition}[section]
\newtheorem{prop}{Proposition}[section]

\title{Billiards in a general domain with random reflections}
\author{Francis~Comets\thanks{Partially supported by CNRS (UMR 7599
``Probabilit{\'e}s et Mod{\`e}les Al{\'e}atoires'')}$^{~,1}$ \and
 Serguei~Popov\thanks{Partially supported by CNPq (302981/02--0),
 and 
``Rede Mate\-m\'atica Brasil-Fran\c{c}a"}$^{~,2}$
\and Gunter~M.~Sch\"utz\thanks{Partially supported by
DFG (Schu 827/5-2, Priority programme SPP 1155)}$^{~,3}$ 
\and Marina Vachkovskaia\thanks{Partially supported 
         by CNPq (306029/03--0 and 200460/06--4)}$^{~,4}$}

\begin{document}

\maketitle

{\footnotesize
\noindent $^{~1}$Universit{\'e} Paris 7, UFR de Math{\'e}matiques,
case 7012, 2, place Jussieu, F--75251 Paris Cedex 05, France\\
\noindent e-mail: \texttt{comets@math.jussieu.fr},
\noindent url: \texttt{http://www.proba.jussieu.fr/$\sim$comets}

\smallskip
\noindent $^{~2}$Instituto de Matem{\'a}tica e Estat{\'\i}stica,
Universidade de S{\~a}o Paulo, rua do Mat{\~a}o 1010, CEP 05508--090,
S{\~a}o Paulo SP, Brasil\\
\noindent e-mail: \texttt{popov@ime.usp.br}, 
\noindent url: \texttt{http://www.ime.usp.br/$\sim$popov}

\smallskip
\noindent $^{~3}$Forschungszentrum J\"ulich GmbH,
Institut f\"ur Festk\"orperforschung,
D--52425 J\"ulich, Deutschland\\
\noindent e-mail: \texttt{G.Schuetz@fz-juelich.de}, \\
\noindent url: \texttt{http://www.fz-juelich.de/iff/staff/Schuetz\_G/}

\smallskip
\noindent $^{~4}$Departamento de Estat\'\i stica, Instituto de Matem\'atica,
Estat\'\i stica e Computa\c{c}\~ao Cien\-t\'\i\-{}fica,
Universidade de Campinas,
Caixa Postal 6065, CEP 13083--970, Campinas SP, Brasil\\
\noindent e-mail: \texttt{marinav@ime.unicamp.br},
\noindent url: \texttt{http://www.ime.unicamp.br/$\sim$marinav}
}

\begin{abstract}
We study stochastic billiards on general tables: a particle moves 
according to its constant velocity inside some domain $\DD \subset \R^d$
until it hits the boundary and bounces randomly inside according
to some reflection law. 
We assume that the boundary of the  domain is locally Lipschitz
and almost everywhere continuously differentiable. The angle of the
outgoing velocity with the inner normal vector
has a specified, absolutely continuous density. We construct 
the discrete time and the continuous time processes recording the 
sequence of hitting points on the boundary and the pair location/velocity.
We mainly focus on the case of bounded domains. Then, we prove exponential
ergodicity of these two Markov processes, we study their invariant 
distribution and their normal (Gaussian) fluctuations. 
Of particular interest is 
the case of the cosine reflection law: the stationary distributions 
for the two processes are uniform in this case, the
discrete time chain is reversible though the continuous time
 process is quasi-reversible. Also in this case, we give a natural 
construction of a chord ``picked at random'' in~$\DD$, and we study
the angle of intersection of the process with a $(d-1)$-dimensional
manifold contained in $\DD$.  
\\[.3cm]{\bf Keywords:} cosine law, Knudsen random walk,
 Knudsen regime, random chord, Bertrand paradox,
kinetic equations, invariant measure, shake-and-bake algorithm
\\[.3cm]{\bf AMS 2000 subject classifications:}
 60J25. Secondary: 37D50, 58F15, 60J10.
\end{abstract}

\section{Introduction}
\label{s_intro}

The purpose of this paper is to provide a rigorous mathematical treatment
of a stochastic process that can be informally described as follows.  A particle moves
with constant speed inside some $d$-dimensional domain. When it hits
the domain boundary, it is reflected in some random direction, not depending on
the incoming direction, and keeping the
absolute value of its speed. One would like to understand basic stationary
and dynamical properties of this process for ``physically reasonable'' domains 
and laws of reflection.

Indeed, the current physics motivation for studying this process comes from 
the need to understand diffusive motion in porous media on a microscopic level.
One wishes to explore the random motion of molecules inside pores
of nanometer scale in order to deduce large-scale transport properties
inside a porous grain. There are numerous applications of this basic
problem. Recently, this has emerged as being of considerable
importance e.g.\ in the design of synthetic zeolites for catalysis \cite{Keil00}
for which in a
computer model a pore-selective molecular traffic control effect has been 
shown  to enhance the effective reactivity of catalytic grains \cite{Brza06}.

Inside a pore only very few molecules can reside so that the mean free
path between molecule-molecule collision is rather large. As a result
intraporous mobility is dominated by the interaction of molecules with the 
pore walls rather than among themselves. This is the so-called Knudsen 
regime in gas dynamics which gives rise to the model described above
in which all molecule-molecule interactions are neglected \cite{Knud52,Cerc88}.
As a further simplifying feature the kinetic energy of a particle is strictly
conserved. Due to a complicated microscopic structure of the pore walls
momentum is assumed to be transferred in a memoryless random fashion. 

Already Martin Knudsen proposed the cosine reflection law, in which case the 
reflected direction has a rotation-invariant cosine distribution around the 
surface normal. 
This is in some sense a physically natural choice
(Feres and Yablonsky \cite{FY,F})
and has attracted very considerable interest in recent years, see e.g.\
Coppens and collaborators  \cite{CD,CM,MC}, and also \cite{RZBK}. 
In particular, it has been debated whether or not 
the self-diffusion coefficient (defined by the mean-square displacement
of a particle) and the transport diffusion coefficient (defined via the
density gradient between two open boundaries of a pore) coincide in the
case of the Knudsen cosine law and how these diffusivities depend on the
structure of the pore walls. In two space dimensions the large scale
motion of a particle inside a strip of constant width is analytically
argued 
to be superdiffusive.

The results reported in \cite{CD,CM,MC,RZBK} on diffusion in more complicated pore 
surfaces are mostly numerical simulation results for open pores where particles are 
injected and extracted. This numerical evidence, however, has proved to be subtle
and often ambiguous, much depending on details of the definition of the
quantities under consideration, which are frequently not provided. In order
to clarify some of these open issues and for making further progress we provide
here a mathematically
rigorous treatment. In the present paper we focus on bounded domains which also prepares the
ground for addressing the problems arising in future study of open domains.
Pores in real solids such as zeolites have a rather complicated surface 
structure which is usually not known in great detail. Hence we keep our
discussion as general as possible as far as the domain boundary and the
reflection law is concerned. Nevertheless, particular emphasis is given to
the previously studied Knudsen cosine reflection law.

Next, we discuss related papers from the probabilistic literature.
Motivated by an asymptotically minimax strategy for the so-called
``princess and monster'' game, stochastic billiards with
the cosine law were considered 
by Lalley and Robbins \cite{LR,LR2}
on a convex domain of the plane, and the invariant measure
for the position was found to be uniform. Later on, billiards 
have been used in Monte Carlo Markov chains algorithms
by Borovkov \cite{B1,B2} and Romeijn~\cite{R} (where the stochastic billiard
with the cosine law of reflection appears as ``running shake-and-bake algorithm'')
to generate random vectors in a connected domain with smooth boundary
or on the boundary itself. S.~Evans~\cite{E} 
gives a detailed account of the case
of uniform reflection law and of a domain with ${\mathcal C}^1$ boundary
or a polygon.
We emphasize that the present paper improves much on the geometry
of the domain.

  From an analytical perspective, Bardos, Golse
and Colonna  \cite{BGC} and Boatto and Golse \cite{BG} consider 
a gas of particles moving along straight lines 
between two infinite parallel horizontal planes, where the reflection
angle is given by  some dynamical system. They prove a diffusion 
limit in this model where there exists a natural angular cut-off. 
For some other models with deterministic behaviour inside the domain and
stochastic behaviour on the boundaries, see Goldstein, Kipnis and Ianiro~\cite{GKI}
and Caprino and Pulvirenti~\cite{CP}.

We now briefly describe the content of this paper.
It turns out to be useful to distinguish the random motion of the particle 
(below 
referred to as stochastic billiard) from the stochastic process defined by 
the sequence
of hits of the particle at the boundary (below referred to as random walk).
We consider domains which boundary is almost
everywhere continuously differentiable and satisfies a Lipschitz condition. 
We first construct the walk and the billiard on such general domains.
We then focus on bounded domains.
For general reflection laws,
we prove for both processes existence and uniqueness of an invariant 
probability measure,
convergence to it, and a central limit theorem. 
A key observation is that these Markov processes satisfy the D\"oblin 
condition. 

The special case of  the cosine law for the  reflection, presents some 
interesting probabilistic, dynamic and geometric properties. The 
invariant measures are uniform, the random walk is reversible,
as well as the stochastic billiard up to a sign change in the speed.
The cosine law also yields a construction 
of a ``chord picked at random'' in the domain, which has
universality properties independently of the domain.
The 
mean chord length is equal, up to an universal constant factor,
to the volume--to--surface ratio  of the domain. 
Also, for a smooth convex subdomain,
our random chord induces a random chord on the subdomain, simply by 
its trace (conditionally on a non-empty intersection): we prove that 
the law of this trace coincides with 
our construction of a random chord when directly 
performed on the subdomain. 
With a view on future applications
to open systems we finally establish some statistical property of the 
particle when it crosses a surface inside the domain.

This paper is organized as follows. In Section~\ref{s_formdefres} we define
the model, comment on conditions for ``good'' behaviour and state our main results. 
Then, in Section~\ref{s_aux} we establish some auxiliary results related mainly
to geometric properties of domains with Lipschitz boundary. Our main results are
proved in Section~\ref{s_proofs_KRW} (discrete time), 
Section~\ref{s_proofs_KSB}
(continuous time) and Section~\ref{s_proofs_RC}
(geometric properties of the random chord).

\section{Formal definitions and main results}
\label{s_formdefres}
In Section~\ref{s_gen_notation} we describe the
general conditions imposed on the domain~$\DD$ and we introduce some notation
related to this domain. In Section~\ref{s_def_model} the discrete time and continuous time
processes are formally defined, and we discuss conditions under which the processes
exist or have good behavior. The results are formulated in Section~\ref{s_mainres}.

\subsection{General notations and standing assumptions}
\label{s_gen_notation}
Let $\|\cdot\|$ be the Euclidean norm in~$\R^d$, and define 
$\BB(x,\eps)=\{y\in\R^d:\|x-y\|<\eps\}$ to be the open $\eps$-neighborhood
of $x\in\R^d$, and $\Sph^{d-1}=\{y\in\R^d:\|y\|=1\}$ to be the unit sphere.
Consider an open connected domain $\DD\subset\R^d$, and let~$\fr$ be 
the boundary of~$\DD$, and~$\cDD$ be the closure of~$\DD$ (so that $\DD\cup\fr=\cDD$). 
In this paper we will deal with three reference measures:
\begin{itemize}
\item the $d$-dimensional Lebesgue measure on $\DD$ or $\bar \DD$,
\item the Haar measure on the sphere $\Sph^{d-1}$,
\item the $(d-1)$-dimensional Hausdorff measure in $\R^d$ restricted to $\fr$.
\end{itemize}
To simplify notations we will use the same 
symbols $dx, dv, dz, \ldots$ to denote all of them. 
We warn the reader, and we recall that ambiguity on
the measures under consideration -- should it arise -- can be easily 
resolved by checking the space of integration.
Similarly, when we write $|A|$, this corresponds to the
$d$-dimensional Lebesgue measure of~$A$ in the case $A\subset\DD$,
and to the $(d-1)$-dimensional Hausdorff measure of~$A$ in the case $A\subset\fr$.

Throughout this paper we suppose that~$\fr$ is a
$(d-1)$-dimensional surface satisfying the Lipschitz condition.
This means that for any~$x\in\fr$ there exist~$\eps_x>0$,
an affine isometry ${\mathfrak I}_x : \R^d\to\R^d$, a function $f_x:\R^{d-1}\to\R$
such that
\begin{itemize}
\item $f_x$ satisfies Lipschitz condition, i.e., there exists a constant~$L_x>0$
such that $|f_x(z)-f_x(z')| < L_x\|z-z'\|$ for all $z,z'$ (without restriction
of generality we suppose that $L_x>1$);
\item ${\mathfrak I}_x x = 0$, $f_x(0)=0$, and
\[
 {\mathfrak I}_x(\DD\cap\BB(x,\eps_x)) = \{z\in\BB(0,\eps_x) : 
                           z^{(d)} > f_x(z^{(1)},\ldots,z^{(d-1)})\}.
\]
\end{itemize}
Recall that, by Rademacher's theorem (cf.\ e.g.~\cite[theorem 3.1.6]{Fed}), 
the Lipschitz condition implies that the boundary~$\fr$ is 
a.e.\ differentiable. This, however, is not enough for our purposes, so we assume additionally
that the boundary is a.e.\ \emph{continuously} differentiable, and we 
denote by $\tRR\subset\fr$ the set of boundary points
where~$\fr$ is continuously differentiable. That is, we suppose that~$\tRR$ is
open (with respect to the induced topology on~$\fr$), such that the $(d-1)$-dimensional
Hausdorff measure of $\fr\setminus\tRR$ is equal to zero, and~$\fr$ has locally~$\cCC^1$
parametrization in any point of~$\tRR$. 

For all~$x\in\tRR$ we can  define a unique
vector $\n(x)\in\Sph^{d-1}$ with the following properties:
\begin{itemize}
\item[(i)] for all small enough $\eps>0$ we have that $x+\eps\n(x)\in\DD$,
\item[(ii)]
\[
 \inf_{y\in\BB(x,\delta) \cap \fr} (y-x)\cdot \n(x) = o(\delta) \qquad \text{ as }\delta\to 0
\]
\end{itemize}
(i.e., $\n(x)$ is the normal vector which points inside the domain). 
Clearly,  the map $x\to\n(x)$ is
continuous for all $x\in\tRR$.

Let $x\cdot y$ be the scalar product of~$x$ and~$y$.
For any $u,v\in\Sph^{d-1}$ we define by
\[
\phi(u,v) = \arccos (u\cdot v)
\]
the angle between~$u$ and~$v$, and by 
\[
\Sph_v = \{w\in\Sph^{d-1}: v\cdot w > 0\}
\]
the half-sphere looking in the direction~$v$.
If $u\in \Sph_{\n(x)}$, we use a simplified notation $\phi_x(u):=\phi(u,\n(x))\in [0,\pi/2)$.

\begin{figure}
\centering
\includegraphics{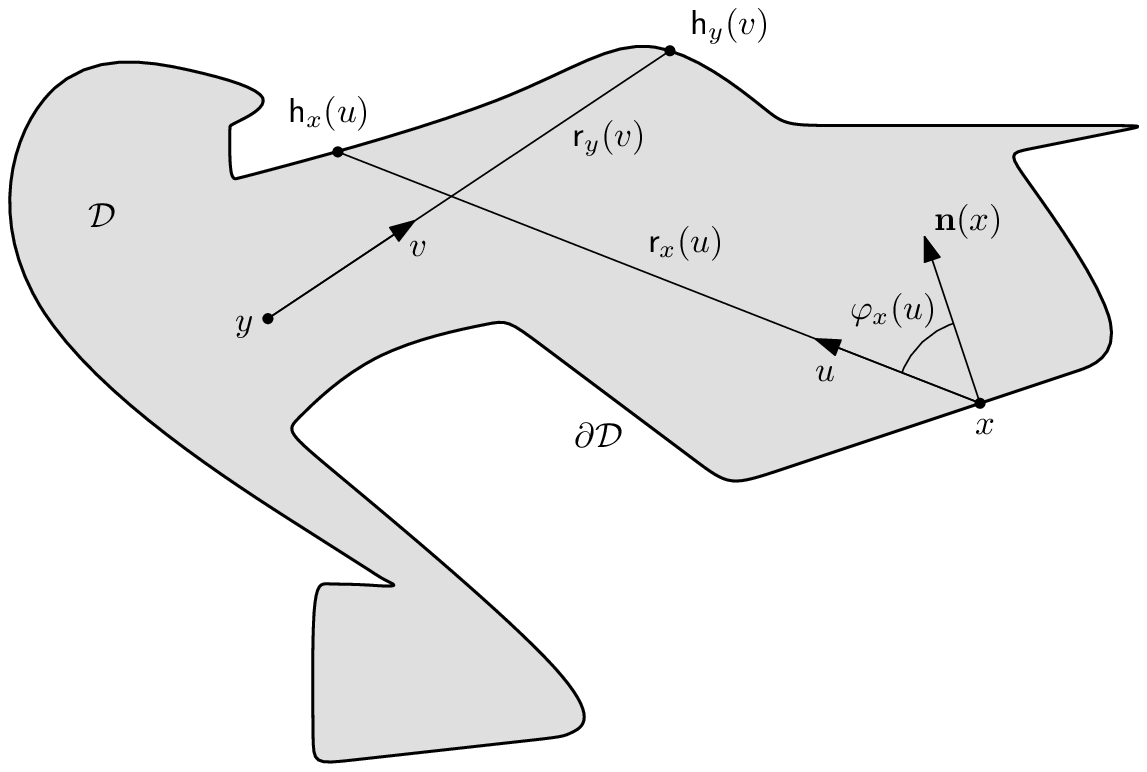}
\caption{On the definition of the quantities $\phi_x(u)$, ${\mathsf r}_x(u)$ 
             and ${\mathsf h}_x(u)$}
\label{def_coisas}
\end{figure}

Note that, by definition of~$\n(x)$, if $u\in\Sph_{\n(x)}$ then for all 
small enough~$t>0$ it holds that $x+tu\in\DD$.
Using this fact, for any pair $(x,u)$ such that either
$x\in\tRR$, $u\in \Sph_{\n(x)}$, or $x\in\DD$, $u\in\Sph^{d-1}$, 
define by (see Figure~\ref{def_coisas})
\[
{\mathsf r}_x(u) = \inf\{t>0 : x+tu \in \fr\}>0
\]
the distance from~$x$ to~$\fr$ along the direction~$u$
(we use the convention $\inf\emptyset=+\infty$), and by
\[
{\mathsf h}_x(u) = x+u{\mathsf r}_x(u) \in \{\fr,\infty\}
\]
the point on~$\fr$ seen from~$x$ in that direction
(in fact, ${\mathsf h}_x(u)=\infty$ if and only if ${\mathsf r}_x(u)=+\infty$).

\subsection{Construction of the walk and billiard}
\label{s_def_model}


In this section we define the main objects of interest, the
continuous time random motion on~$\cDD$ described informally above,
and its skeleton random walk, i.e., the discrete-time
 Markov chain corresponding to the sequence of points where the 
continuous time process hits~$\fr$.

Denote by $e=(1,0,\ldots,0)$ the first coordinate vector.
Let $\eta=(\eta_n,n=1,2,\ldots)$ be a sequence of i.i.d., $\Sph_{e}$-valued
random vectors with absolutely
continuous distribution on~$\Sph_{e}$:
there is a function~$\bga$ such that
for any measurable $B\subset\Sph_{e}$ it holds that 
$\PP[\eta\in B]=\int_{B}\bga(u) \, du$ (where~$\PP$ stands for the probability). 
We suppose that~$\bga$ is supported
on~$\Sph_{e}$ in the  sense that
\begin{equation}
\label{support_halfsph}
 \inf_{u\in B}\bga(u) > 0
\end{equation}
for any compact $B\subset\Sph_{e}$.

To define the model, we need also a family $(U_x,x\in\fr)$ of unitary linear operators in~$\R^d$
with the property $U_x e = \n(x)$ for all $x\in\tRR$, 
so that~$U_x$ is a rotation of~$\Sph^{d-1}$
which sends~$e$ to~$\n(x)$. Here, two cases need to be considered separately.
\begin{itemize}
\item Case~1: the function~$\bga$ is rotationally symmetric around the 
first coordinate axis, i.e, $\bga(u_1)=\bga(u_2)$
 whenever $u_1\cdot e = u_2 \cdot e$. In this case we 
do not need any additional assumptions
 on the family $(U_x,x\in\fr)$ except measurability 
of $x \mapsto U_x$, and we denote $\RR:=\tRR$.
\item Case~2: the function~$\bga$ is not symmetric. In this case we assume also that
the map $x\to U_x$ is continuous for almost all $x\in\fr$, and we define
$\RR=\{x\in\tRR : \text{ the map } x\to U_x \text{ is continuous in }x\}$.
\end{itemize}

Fix an arbitrary $x_0\in\RR$; using the sequence $(\eta_n,n\geq 1)$, we define
a discrete time Markov chain $(\xi_n,n=0,1,2,\ldots)$ and a continuous time
process $(X_t,t\geq 0)$, so-called  \emph{stochastic billiard},
in the following way.
 Below, ${\mathfrak S}$ stands for a special state, that we
allow for technical reasons, 
$\xi_n$ or $X_t={\mathfrak S}$ means that the corresponding
process is undefined.
\begin{itemize}
\item[(i)] Put $\xi_0=x_0$.
Inductively, we define $\xi_{n+1}\in\{\fr,\infty,{\mathfrak S}\}$ as follows:
\[
 \xi_{n+1} = \left\{\begin{array}{ll}
                  {\mathsf h}_{\xi_n}(U_{\xi_n}\eta_n), & \text{if }\xi_n\in\RR,\\
                  {\mathfrak S}, &\text{if }\xi_n={\mathfrak S} 
                       \text{ or }\xi_n\in\fr\setminus\RR,\\
                  \infty, &\text{if }\xi_n=\infty.
                \end{array}\right.
\]
Since $\eta_n$ has a density, one can show 
(see Lemma~\ref{l_abs_continuity} below) that for any $x\in\RR$, 
the distribution of~${\mathsf h}_{x}(U_{x}\eta)$
is absolutely continuous with respect to the $(d-1)$-dimensional Hausdorff measure on~$\fr$.
Since almost all points of~$\fr$ belong to~$\RR$, almost surely it holds
that either $\xi_n\in\fr$ for all~$n$, or there exists some~$n_0$ such that
$\xi_{n_0}=\infty$.
\item[(ii)] Again, we put $X_0=x_0$, and let us abbreviate $\tau_0:=0$ and
\[
 \tau_n = \sum_{i=0}^{n-1} {\mathsf r}_{\xi_i}(U_{\xi_i}\eta_i),
\]
$n=1,2,\ldots$ (formally, if $\xi_k={\mathfrak S}$ for some~$k$, we put $\tau_{k'}=\tau_k$
for all $k'\geq k$). 
Then, for~$t\in [\tau_n,\tau_{n+1})$ define
\[
 X_t = \xi_n + (t-\tau_n)U_{\xi_n}\eta_n.
\]
Since a.s.\ the random walk never enters $\fr\setminus\RR$, we obtain that, with probability~$1$,
$X_t\in\cDD$ is well defined for all~$t\in [0,\lim_{n\to\infty}\tau_n)$.
\end{itemize}
Intuitively, the Markov chain~$\xi$ 
can be described as follows: on each step, the particle chooses
the direction of the jump according to the density~$\bga$ ``turned'' in the
direction of the normal vector (the one that looks inside the domain~$\DD$),
and then instantly jumps to the point on the boundary of~$\DD$ seen in this direction. 
As a process on the boundary of $\DD$, it is as close as possible 
to a random walk, and we will denote it like that.
For some unbounded domains it may happen, of course, 
that such a point does not exist, 
in which case the particle jumps directly to the infinity. For 
the stochastic billiard~$X$, when on the
boundary, the particle uses the same procedure to choose the direction, but then it moves
with speed~$1$ until it hits the boundary again (if it ever does).

\begin{figure}
\centering
\includegraphics{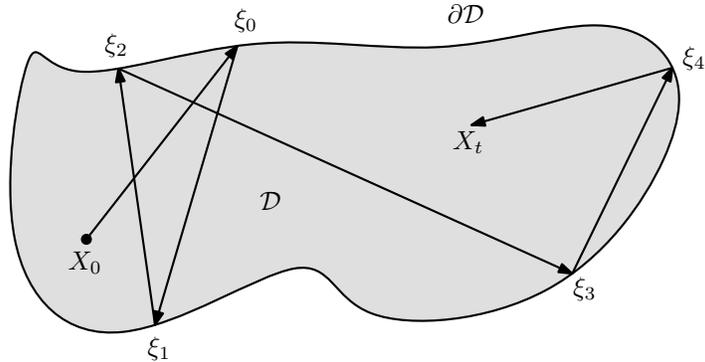}
\caption{A stochastic billiard starting from the interior of~$\DD$}
\label{fig_billiard}
\end{figure}
Let us note that we can define the stochastic billiard starting in any $x_1\in\DD$ as well
(in this case we have to specify also the initial direction), see Figure~\ref{fig_billiard}. Of course,
one may consider the random walk~$\xi$
starting from any initial distribution on~$\fr$ 
(one has to require that the probability that the starting point belongs
to~$\RR$ is~$1$); similarly, one may consider the stochastic billiard~$X$ starting from
any initial distribution on $\DD\times\Sph^{d-1}$
(and, in fact, we will frequently do that in this paper).

The next proposition shows that, supposing that the boundary of the 
domain satisfies the Lipschitz condition,
the process~$X_t$ is well defined for all~$t\geq 0$.
\begin{prop}
\label{tau_to_infty}
With probability~$1$, $\lim_{n\to\infty}\tau_n = +\infty$.
\end{prop}

It would be interesting to consider more general laws of reflection: 
for example,
one can suppose that the law of~$\eta$ is continuous, but is not supported on
the whole~$\Sph_{e}$; and/or one can suppose that the law of~$\eta$ has discrete
and/or singular components. However, this would bring in some substantial technical
difficulties. For instance, in the former case the domain may be such that some
parts of the boundary are not seen from any of the boundary points (since we no longer
have the rule ``if the particle can go from~$x$ to~$y$, 
then it can go from~$y$ to~$x$''). In the latter case, there
is an additional difficulty that, even if~$x\in\RR$, it may happen that 
$\PP[{\mathsf h}_{x}(U_{\n(x)}\eta)\in\RR]<1$, and so one needs to define the 
set~$\RR$ differently (to ensure that with probability~$1$ 
the random walk never leaves the set~$\RR$), 
and show that, with that new definition, we still have $|\fr\setminus\RR|=0$.

\begin{figure}
\centering
\includegraphics{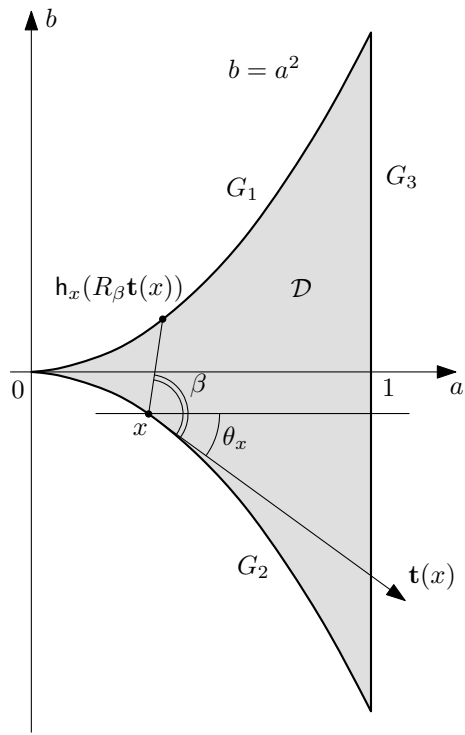}
\caption{A stochastic billiard on a non-Lipschitz domain}
\label{fig_x2}
\end{figure}

We end this subsection by commenting our assumption
that the boundary satisfies
Lipschitz condition. If one considers domain with the boundary which is
not everywhere continuously differentiable (even if one supposes that it is
not differentiable in only one point), then the D\"oblin condition~(\ref{Doeblin})
below need not be satisfied when the boundary is non-Lipschitz (for instance, it is not
satisfied for the domain in Figure~\ref{fig_x2}). So, although it can be
indeed interesting to study stochastic billiards on such domains, this task
is much more difficult due to the lack of geometric ergodicity and related
properties (in~\cite{B1} it was claimed that the results there are valid 
for domains with boundary consisting of a finite number of smooth surfaces;
on the other hand, Theorem~1 of~\cite{B1} is wrong for e.g.\ the domain 
of Figure~\ref{fig_x2}). Also, without Lipschitz condition it is not even
clear if the process~$X_t$ is well defined for all~$t\geq 0$. If the density~$\bga$
is not symmetric, there are examples (cf.\ the one below) when the random walk~$\xi_n$
converges to a (singular) point on the boundary. It is, in our opinion, a
challenging problem to find out if such examples exist for symmetric
laws of reflection.

\medskip
\noindent
{\bf Example (Non-Lipschitz boundary).} Consider the domain 
\[
\DD=\{(a,b)\in\R^2 : 0<a<1, -a^2<b<a^2\}
\]
(see Figure~\ref{fig_x2}). We denote by
\begin{eqnarray*}
 G_1 &=& \{(a,b)\in\R^2 : 0\leq a<1, b=a^2\},\\
 G_2 &=& \{(a,b)\in\R^2 : 0<a<1, b=-a^2\},\\
 G_3 &=& \{(a,b)\in\R^2 : a=1, 0\leq b \leq 1\},
\end{eqnarray*}
the three smooth pieces of the boundary, 
so that $\fr=G_1\cup G_2\cup G_3$.
Let $x=(x_1,-x_1^2)\in G_2$ be a point on the lower part of the boundary,
and abbreviate $\theta_x=\arctan (2x_1)$ (we use the same notation for
$x=(x_1,x_1^2)\in G_1$ as well). Let~$R_\beta$ be the operator
of anticlockwise rotation by angle~$\beta \in [0,\pi]$, 
and~$J$ be the symmetry
operator defined by $J(a,b)=(-a,b)$. The family of operators~$(U_x,x\in\fr)$
is then defined as follows: $U_x=R_{\frac{\pi}{2}+\theta_x}J$ for $x\in G_1$,
$U_x=R_{\frac{\pi}{2}-\theta_x}$ for $x\in G_2$, 
$U_x=J$ for $x\in G_3$.

Let $g(x)$ be the first coordinate of $x\in\R^2$. 
Define, for $x \in G_1 \bigcup G_2$,
\[
 D_x(\beta) = g({\mathsf h}_x(R_\beta{\mathbf t}(x)))-g(x),
\]
where~${\mathbf t}(x)$ is the tangent unit vector at~$x$ with positive
first coordinate.

After some elementary computations, one can obtain the following for small enough~$x$:
\begin{itemize}
\item $D_x(\beta)\leq 1$ for $\beta < 8x_1$,
\item $D_x(\beta)\leq \frac{8x_1^2}{\beta}$ for $\beta\in [8x_1;\frac{\pi}{2}+\theta_x]$,
\item $D_x(\beta)\leq 0$ for $\beta\in [\frac{\pi}{2}+\theta_x;\frac{3\pi}{4}+\theta_x]$,
\item $D_x(\beta)\leq -\frac{2x_1^2}{3}$ for $\beta >\frac{3\pi}{4}+\theta_x$.
\end{itemize}

Consider the following density on~$\Sph_e$: being~$\beta$ the
angle between the element of~$\Sph_e$ and $(-e_2)$, write 
(with some abuse of notation) $\bga(\beta)=\frac{4}{\pi^4}\beta^3$,
$\beta\in (0,\pi)$. 
Since
\[
 \E(g(\xi_{n+1})-g(\xi_n)\mid \xi_n=x) = \frac{4}{\pi^4}\int_0^\pi D_x(\beta)\beta^3\,d\beta,
\]
one sees from the above above properties of~$D_x(\beta)$ that
\begin{equation}
\label{supermart}
 \E(g(\xi_{n+1})-g(\xi_n)\mid \xi_n=x) \leq 0
\end{equation}
for all small enough~$x$. Define $\sigma(y_0)=\min\{n:g(\xi_n)\geq y_0\}$.
 The equation~(\ref{supermart}) means that, for small enough $y_0$,
the process $g(\xi_{n\wedge\sigma(y_0)})$ is a supermartingale. So, there
exists $g_\infty = \lim_{n\to\infty}g(\xi_{n\wedge\sigma(y_0)})$, and,
by the Fatou lemma we have $\E g_\infty \leq \E g(\xi_0)$. If~$\xi_0<y_0$,
this implies that $\sigma(y_0)=\infty$ with positive probability, and then with
a little more work one concludes that 
\[
\xi_n\longrightarrow 0 \; {\rm a.s.\ \ as\ } n \to \infty.
\]
It seems clear that one should be able to construct examples of such kind
where it holds also that
 $\lim_{n\to\infty}\tau_n<\infty$. However in this particular example
it is rather difficult to find out whether this happens, so we
opted for not discussing it here. Anyway, existence of such examples is 
already a good reason to assume that the boundary of the domain satisfies 
the Lipschitz condition.
\qed

\subsection{Main results}
\label{s_mainres}

For the rest of this paper, we suppose that~$\DD$ is a bounded domain 
(that is, $\diam(\DD)<\infty$).
Note that the Lipschitz condition on the boundary implies that 
in this case it holds that $|\fr|<\infty$
and that the number of connected components of~$\fr$ is finite, see Lemma~\ref{l_basic_Lip} below.

For any $x,y\in\R^d$, $x\neq y$, let us define the vector
$\ell_{x,y}=\frac{y-x}{\|y-x\|}\in\Sph^{d-1}$; note that $\ell_{x,y}=-\ell_{y,x}$.
We say that $y\in\cDD$ is {\it seen from} $x\in\cDD$ if
there exists $u\in\Sph^{d-1}$ and~$t_0>0$ such that $x+tu\in\DD$ for all $t\in (0,t_0)$
and $x+t_0 u = y$. Equivalently, the open segment $(x,y)$ lies in~$\DD$ and 
$x\neq y$.
Clearly, if~$y$ is seen from~$x$ then~$x$ is seen from~$y$,
and we write ``$x \leftrightarrow y$'' when this occurs. 

\begin{figure}
\centering
\includegraphics{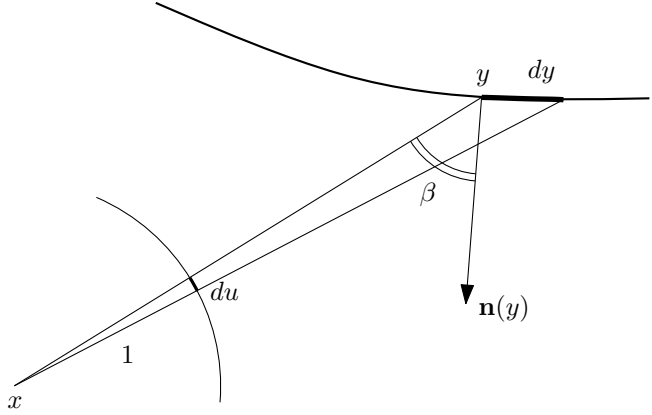}
\caption{$du=\|x-y\|^{d-1}\cos \beta\, dy$}
\label{fig_differentials}
\end{figure}

   From the definition of the Markov chain~$\xi_n$,  we will show that 
the transition kernel can be written
\[
\PP[\xi_{n+1}\in A \mid \xi_n=x]=\int_A K(x,y)\, dy\;,
\]
where (see Figure~\ref{fig_differentials})
\begin{equation}
\label{def_trans_dens}
 K(x,y) = \frac{\bga(U_x\ell_{x,y})\cos\phi_y(\ell_{y,x})}{\|x-y\|^{d-1}}
\1{x,y\in\RR, x \leftrightarrow y}
\end{equation}
is the transition density.
Formally, if~$x$ or~$y$ does not belong to~$\RR$, we put $K(x,y)=0$.

As mentioned above, we consider separately the important case of the cosine law of reflection
\begin{equation}
\label{cosine_density}
 \bga(u) = \gamma_d\cos\phi(e,u) = \gamma_d e\cdot u,
\end{equation}
where the normalization constant~$\gamma_d$ is defined by
\[
\gamma_d = \Big(\int_{\Sph_{e}} \cos\phi(e,u)\, du\Big)^{-1}
         = \Big(\int_{\Sph_{e}} u\cdot e\, du\Big)^{-1}\;,
\]
making $K(x,\cdot)$ a probability density on $\fr$ for all $x \in \RR$.
The transition density then becomes
\begin{eqnarray*}
 \tK(x,y) &=& \frac{\gamma_d\cos\phi_x(\ell_{x,y})\cos\phi_y(\ell_{y,x})}{\|x-y\|^{d-1}}
\1{x,y\in\RR, x \leftrightarrow y}\\
 &=& \frac{\gamma_d \big((y-x)\cdot\n(x)\big)\big((x-y)\cdot\n(y)\big)}{\|x-y\|^{d+1}}
\1{x,y\in\RR, x \leftrightarrow y}
\end{eqnarray*}
When considering the specific case of the cosine law of reflection, we write $\txi,\tX$
instead of $\xi,X$, and we call these processes \emph{Knudsen random walk} (KRW)
and \emph{Knudsen stochastic billiard} (KSB).

Define~${\hat \mu}_0$
to be the ``uniform'' probability
measure on~$\fr$: ${\hat \mu}_0(A) = \frac{|A|}{|\fr|}$.
Since $\tK(x,y)$ is symmetric, we immediately
obtain that the KRW~$\txi$ is reversible,
with the reversible (and thus invariant) measure~${\hat \mu}_0$.
For other reflection laws it is usually not easy to find the exact form of the
invariant measure (except for some particular cases, see~\cite{E}), but nevertheless
we prove that such a measure exists and is unique, the random walk converges to it
exponentially fast, and the Central Limit Theorem holds.
\begin{theo}
\label{t_conv_KRW}
Suppose that $\diam(\DD)<\infty$.
\begin{itemize}
\item[(i)] There exists a unique
probability measure~${\hat \mu}$ on~$\fr$ which is invariant for the
random walk~$\xi_n$. Moreover, there exists a function~$\psi: \fr\to\R_+$
such that ${\hat \mu}(A) = \int_A\psi(x)\,dx$
which satisfies
\begin{equation}
\label{eq_stat_dens}
 \psi(x) = \int_{\fr}\psi(y)K(y,x)\,dy.
\end{equation}
Finally, the density $\psi$ can be chosen in such a way that $\inf_{\fr}\psi>0$.
\item[(ii)] There exist positive constants~$\beta_0,\beta_1$
(not depending on the initial distribution of $\xi_0$) such that
\begin{equation}
\label{eq_conv_KRW}
 \|\PP[\xi_n\in \cdot] -{\hat \mu}\|_{\sf v} \leq \beta_0 e^{-\beta_1 n},
\end{equation}
where $\|\cdot\|_{\sf v}$ is the total variation norm.
\item[(iii)] We have the Central Limit Theorem:
for any measurable $A\subset\fr$ there exists $\sigma_A$, with
$\sigma_A> 0$ if $0<|A|<|\fr|$, such that
\begin{equation}
\label{eq_CLT_KRW}
 n^{-1/2} \Big(\sum_{i=1}^n \1{\xi_i\in A} - n{\hat \mu}(A)\Big)
\end{equation}
converges in distribution to a $\text{\rm Normal}(0,\sigma_A^2)$ random variable.
\end{itemize}
In particular, for the KRW, (\ref{eq_conv_KRW}) and~(\ref{eq_CLT_KRW}) hold
with ${\hat \mu}={\hat \mu}_0$.
\end{theo}
We remark that the relation (\ref{eq_stat_dens}) for the Radon-Nykodim density~$\psi$
of the invariant measure~${\hat \mu}$
is what is usually called the balance condition.

Contrary to~$\xi_n$, the process~$X_t$ by itself is not Markovian,
which, in principle, makes it more difficult to analyze. To overcome that difficulty,
we define another process~$V_t$ by
\[
 V_t= \lim_{\eps\downarrow 0}\frac{X_{t+\eps}-X_{t}}{\eps}.
\]
Clearly, $V_{t_1}=V_{t_2}=U_{\xi_{n-1}}\eta_{n-1}$
for all $t_1,t_2\in [\tau_{n-1},\tau_n)$. In words, $V_t$ is
the c\`adl\`ag version of
the motion direction (or  speed vector) of the stochastic billiard process
 at time~$t$. As before, for the particular
case of the cosine reflection law, we write $\tV$ instead of~$V$.

Now, it is clear that the pair $(X_t,V_t)$ is a Markov process.
Associated to this Markov process there is a whole family of martingales.
Introduce the sigma-field $\FF_t$  generated by $X_s,V_s, s\in[0,t]$,
and $ \LL $ the transport operator
\[
 \LL f(x,v) = v\cdot \nabla_x f.
\]
Let $\cCC^{1,0}(\bar \DD \times \Sph^{d-1})$ be the set of functions
on $\bar \DD \times \Sph^{d-1}$ which are continuous and once  continuously
differentiable  in the first variable on $\bar \DD \times \Sph^{d-1}$.
\begin{prop}
\label{prop:pbmg}
Let $f \in \cCC^{1,0}(\bar \DD \times \Sph^{d-1})$
such that
\begin{equation} \label{eq:dom}
f(x,v)=\int_{\Sph_e} f(x, U_{x}u)  \bga(u) du \qquad
\forall x \in \RR, v\cdot \n(x)\leq 0.
\end{equation}
Then,
\begin{equation} \label{eq:pbmart}
f(X_t,V_t)- \int_0^t  \LL f(X_s,V_s)ds
\end{equation}
is a martingale with respect to $(\FF_t; t\geq 0)$.
\end{prop}

 From this, it follows that the  generator of the  Markov process
$(X_t,V_t; t\geq 0)$ is
given by~$\LL$ on the set of smooth functions satisfying~(\ref{eq:dom}).
Note that the left-hand side of~(\ref{eq:dom}) does not depend on~$v$ 
provided $x \in \RR, v\cdot \n(x)\leq 0$.

\medskip
Next, we show that  the stochastic billiard converges exponentially fast
to equilibrium.
\begin{theo}
\label{t_conv_KSB}
There exist a probability measure~$\chi$ on
$\DD\times\Sph^{d-1}$ and positive constants~$\beta'_0,\beta'_1$ 
(not depending on the initial distribution of position and direction) such that
\begin{equation}
\label{eq_conv_KSB}
 \|\PP[X_t\in \cdot, V_t\in \cdot] - \chi\|_{\sf v} \leq \beta'_0 e^{-\beta'_1 t},
\end{equation}
for all~$t\geq 0$.
\end{theo}

\medskip

We complement this result with a central limit theorem.
\begin{theo}
\label{t_cl_KSB}
Let $f: \bar \DD \times \Sph^{d-1} \to \R$ be measurable and bounded function.
Suppose also that it is centered, that is,
\begin{equation}
\label{eq:centered}
 \int_{ \bar \DD \times \Sph^{d-1}} f(x,v) \chi(dx,dv) =0.
\end{equation}
Then, as $t \to \infty$, the random variables 
\[
 t^{-1/2} \int_0^t f(X_s,V_s)\, ds
\]
converge in distribution to a $\text{\rm Normal}(0,\sigma_f^2)$
random variable, $\sigma_f^2 \geq 0$.
\end{theo}

\noindent
{\bf Remark.} The variance $\sigma_f^2$ may be equal to zero
for some non-trivial $f$'s. Here is an example for KSB.
Let $G: \bar \DD \to \R$ be a $\cCC^1$ function, and
\begin{equation}
\label{func_var0}
 f(x,v)= - v \cdot \nabla G(x).
\end{equation}
Then, $\sigma_f^2=0$. The proof of this fact is placed in the end of the
proof of Theorem~\ref{t_cl_KSB}. \qed

\medskip

Next, we relate the invariant measure for the stochastic billiard
to the one for the skeleton random walk. Recall that, by Theorem~\ref{t_conv_KRW},
the measure~$\hat \mu$ is absolutely continuous
with respect to the Hausdorff measure~$dz$ on~$\fr$, and $\psi(z):=\frac{d \hat \mu}{dz}(z)$.
Let $\mu_0$ and $\nu_0$ be the uniform measures on~$\DD$ and~$\Sph^{d-1}$ respectively,
i.e., $\mu_0(A) = \frac{|A|}{|\DD|}$, $\nu_0(B) = \frac{|B|}{|\Sph^{d-1}|}$.
\begin{theo}
\label{th_gen}
The invariant measure~$\chi$ is absolutely continuous with respect
to $\mu_0\otimes\nu_0$, and is given by
\[
\chi(dx,dv)= \psi(z)\;
\frac{\bga( U_z^{-1} v) }{\cos \phi_z(v)}dx\,dv, \qquad z={\mathsf h}_{x}(-v).
\]
In particular,
the product measure $\mu_0\otimes\nu_0$ is invariant for
$KSB$ $(\tX_t,\tV_t)$.
\end{theo}

For the rest of this section, we concentrate on the process
with the cosine reflection law.

By means of an usual procedure (shift the law of the KSB starting with the stationary
measure $\mu_0\otimes\nu_0$ by $(-T)$, and take the limit as $T\to\8$)
we can define a stationary version of the KSB for all $t\in\R$.
\begin{theo}
\label{th_gen2}
Let $(\tX_t,\tV_t)_{t \in \R}$ be the stationary KSB,
i.e.,
$(\tX_t,\tV_t)$ is uniformly distributed on $\DD \times \Sph^{d-1}$.
Then,
\begin{equation}
\label{eq:quasirev}
 (\tX_t,\tV_t)_{t \in \R} \eqlaw (\tX_{-t},-\tV_{-t})_{t \in \R}.
\end{equation}
\end{theo}

This relation shows that KSB with uniform initial condition
is not only stationary, but also
reversible up to a sign change in the speed.
Also, it yields a pathwise construction of the 
stationary version of the KSB for all $t\in\R$
in the following way: first, choose a random point $(x_0,v_0)\in\cDD\times\Sph^{d-1}$
according to the measure $\mu_0\otimes\nu_0$. Construct a KSB $(\tX^{(1)}_t,\tV^{(1)}_t)$
with the initial condition $(\tX^{(1)}_0,\tV^{(1)}_0)=(x_0,v_0)$, and an independent
KSB $(\tX^{(2)}_t,\tV^{(2)}_t)$ with the initial condition $(\tX^{(2)}_0,\tV^{(2)}_0)=(x_0,-v_0)$.
Then, 
the process  $(\tX_t,\tV_t)_{t\in\R}$ defined by
\[
(\tX_t,\tV_t) = 
 \left\{ 
\begin{array}{cc}
 (\tX^{(1)}_t,\tV^{(1)}_t)  & \quad \mbox{for}\; t\geq 0 \\
&\\
(\tX^{(2)}_{-t},-\tV^{(2)}_{-t})   & \quad \mbox{for}\;t\leq 0
\end{array}
\right.
\]
is a  stationary KSB.

\medskip

For any measurable $A\subset\DD$, define
\[
 \m(A) = \E_{{\hat \mu}_0}\int_0^{\tau_1}\1{\txi_0+ t\ell_{\txi_0,\txi_1} \in A}\, dt
\]
so $\m(A)$ is the mean time that KSB
starting from~$\txi_0$ which is placed uniformly at random on~$\fr$,
spends in~$A$ till the next  hit of the boundary.
In particular, $\bm:=\m(\DD)$ is the mean time interval between the consecutive hittings of the
boundary. Use the abbreviation
\[
 \kappa_d:=\gamma_d|\Sph^{d-1}| = \frac{\pi^{1/2}\Gamma(\frac{d+1}{2})d}{\Gamma(\frac{d}{2}+1)}
\]
(for instance, $\kappa_2=\pi$, $\kappa_3=4$).
 Next, we calculate explicitly
the quantity $\m(A)$ (for the case of convex smooth domain, this formula can be obtained
using the results of~\cite{B1}).
\begin{theo}
\label{th_chord}
Suppose that $\diam(\DD)<\infty$. We have
 \begin{equation}
\label{mean_chord}
 \m(A) = \kappa_d\frac{|A|}{|\fr|}.
\end{equation}
In particular, $\bm=\kappa_d\frac{|\DD|}{|\fr|}$.
\end{theo}

The quantity $\bm$ can be interpreted as ``the mean chord length''
for the domain~$\DD$, when the random chord constructed as follows: 
take a point on~$\fr$ uniformly at random,
and draw a line from there using the cosine probability distribution.
Formally:
\begin{df}
\label{def_random_chord}
The random chord for a bounded domain~$\DD$ satisfying our standing 
assumptions of Section \ref{s_gen_notation},
is a pair of random variables $(\Xi_1,\Xi_2)$,
in $\fr$, 
with the joint density $|\fr|^{-1}\tK(x,y)$. 
\end{df}
Note that, by symmetry of the above density, the construction defines 
an unoriented chord, 
as one expects from a random chord. 
\medskip

\noindent
{\bf Remark.} Here is another construction of the random chord,
which is equivalent to Definition~\ref{def_random_chord}:
pick some $(x,v)$ in $\DD \times \Sph^{d-1}$ according to the density
$\frac{C}{\|{\mathsf h}_x(v)-{\mathsf h}_x(-v)\|} dx \, dv$ 
($C$ is the suitable normalizing constant), and put
$\Xi_1:={\mathsf h}_x(v)$, $\Xi_2:={\mathsf h}_x(-v)$.
This is so because the KSB beginning with a uniform distribution on~$\fr$
is a Palm version of the stationary KSB, and therefore the above fact follows from
Theorem~4.1 of Chapter~8 of~\cite{T} and Theorem~\ref{th_gen}. \qed
\medskip

With this definition, 
the quantity $\bm=\E\|\Xi_1-\Xi_2\|$ is 
the mean chord length, and Theorem
\ref{th_chord} has a nice geometric interpretation: 
the area (volume, \dots) of the domain is always proportional to the 
product of its perimeter (surface area, \dots) and~$\bm$:
\begin{equation}
\label{vol=per*meanchord}
 |\DD| = \frac{\bm |\fr|}{\kappa_d}.
\end{equation}
This result also shows that one can obtain the quantity~$\frac{|\DD|}{|\fr|}$
directly (i.e., without obtaining separately~$|\DD|$
and~$|\fr|$) by simulating the KRW~$\txi_n$.

\medskip

Let us discuss some further properties of the above definition of random chord. The term 
``random  chord''
should remind the reader of the so-called Bertrand Paradox 
(see, for example,~\cite{Jaynes}), which can be briefly described as follows. 
A chord is chosen at random in a circle. The question is, what is the 
probability that the chord is longer than a side of the equilateral triangle 
inscribed in this circle?
To answer this question, one considers three natural methods for choosing a chord at random.
\begin{itemize}
\item[1.] The ``random endpoints'' method: choose two points independently uniformly  on the circumference and draw the chord joining these two points. 
Note that this is equivalent to choosing a uniform point on the circumference and drawing a chord
from there with the angle (between the chord and the radius) chosen uniformly
at random in $(-\pi/2,\pi/2)$.
\item[2.] The ``random radius'' method: choose uniformly 
a radius of the circle and a point on the radius and 
construct the chord whose midpoint is the chosen point. It is easily verified that 
this is equivalent to choosing a uniform point on the circumference
and drawing a chord from there according to the cosine distribution, 
so 
this method is equivalent to our method of choosing a random 
chord on a circle.
\item[3.] The ``random midpoint'' method: choose a point uniformly at random within the 
circle and construct a chord with the chosen point as its midpoint.
\end{itemize}
By elementary computations, one finds 
that the above probability is equal to $1/3$, $1/2$, $1/4$ 
if the  random chord is drawn from methods 1,2 and 3 respectively.
It is then 
routinely argued, that the problem is not well formulated, since it is not clear which 
of the methods should be preferred, and the term ``random chord'' itself is not well defined. 
However, E.~Jaynes~\cite{Jaynes} 
 suggested that, among the three above methods, the method 2 is the 
``right'' one, since it has certain invariance
properties.  (Recall that the method 2 is equivalent to our definition 
of random chord for the particular 
case of a circle.)
In the next theorem we show that the random chord in the sense
of Definition~\ref{def_random_chord} also has similar properties, 
thus leading us to suggest that 
perhaps 
Definition~\ref{def_random_chord} is 
the ``right'' definition of a random chord for general domains.

\begin{figure}
\centering
\includegraphics{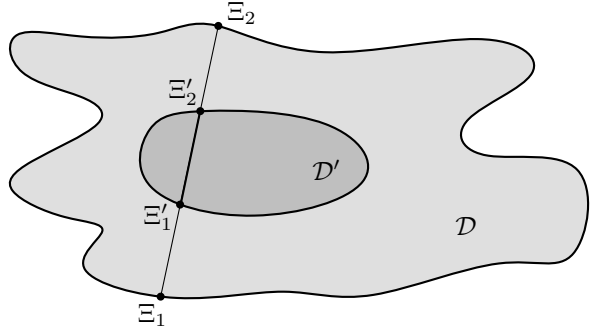}
\caption{A random chord on~$\DD'$ induced by a random chord on~$\DD$}
\label{restr_chord}
\end{figure}

If  $\DD'\subset\DD$ is convex, a  chord 
of~$\DD$ which intersects~$\DD'$ defines a unique chord on~$\DD'$
by its intersection (see Figure~\ref{restr_chord}). 
Let us generate 
independent random chords $(\Xi_1(i),\Xi_2(i))$ of~$\DD$, $i=1,2,\dots$,
till the chord hits the domain~$\DD'$, and then denote by  $(\Xi_1',\Xi_2')$
the intersection. We call  $(\Xi_1',\Xi_2')$ the {\it induced} chord on~$\DD'$.
Like in the ``acceptance-rejection'' algorithm (or ``hit or miss'')
for random variable simulation, 
we easily check that  $(\Xi_1',\Xi_2')$ has the same law as
the endpoints of  $[\Xi_1,\Xi_2]\cap \DD'$ given that these sets intersect.

\begin{theo}
\label{t_restr_convex} 
Let $\DD'\subset\DD$ be convex.
Then, the  chord  induced on~$\DD'$ by the random chord  of~$\DD$
is  the random chord of~$\DD'$ in the sense 
of Definition~\ref{def_random_chord}. 
\end{theo}

When the domain $\DD'$ is not convex, the situation is more complicated, 
since a chord of $\DD$
can intersect $\DD'$ in several places. Denote by 
\[
 (\Xi_{1,k}',\Xi_{2,k}')\;, \quad k \leq \iota
\]
the chords obtained by intersecting~$\DD'$ 
with the chord $(\Xi_1,\Xi_2)$ of $\DD$, with 
$\iota \in \{0,1,\ldots, \infty\}$
their number. 

\begin{theo}
\label{t_nonconvex}
Let $\DD'\subset\DD$ be with the Lipschitz condition and an
almost everywhere continuously differentiable boundary. 
For every Borel subset~$C$  of $\partial\DD'\times \partial\DD'$, we have
\begin{equation} 
\label{eq:baila}
\E \Big[ \sum_{k =1}^{ \iota} \1{(\Xi_{1,k}',\Xi_{2,k}')\in C}
\Big] = \frac{|\fr'|}{|\fr|} \times {|\fr'|}^{-1}
\int_{C}\tilde K(x,y)dx\, dy.
\end{equation}
In particular, 
\begin{equation} 
\label{eq:baila1}
\E [ \iota ] = \frac{|\fr'|}{|\fr|}.
\end{equation}
\end{theo}
Hence, the relation (\ref{eq:baila}) has the following interpretation:
the expected value of the measure giving unit weight to each chord induced
on $\DD'$ by the $\DD$-random chord, 
is the product of the expected  number of induced chords
and the probability distribution for the (ordered) 
endpoints of the random chord of~$\DD'$. 
Another comment on this theorem is
that the number $\iota$ of induced chords is integrable, a property
which does not 
seem easy to prove directly.
The formula~(\ref{eq:baila1}) could remind the reader of the following fact 
concerning the so-called Poisson line process (see e.g.\ Section 8.4.2 of~\cite{SKM}): 
if~$\DD$ is compact and convex, then the number of lines hitting~$\DD$ has 
Poisson distribution with the mean proportional to the perimeter of~$\DD$.


\medskip

Finally, we are interested in what happens at the moments when the KSB crosses a 
surface~$S\subset \DD$.
Suppose that~$S$ is a compact (but not necessarily connected) $(d-1)$-dimensional manifold
(possibly with boundary). Similarly to what we did for~$\fr$,
 we assume also that~$S$ is locally Lipschitz
and almost everywhere continuously differentiable. 
So, for a.a.\ $x\in S$ we can define 
a normal vector~$\n(x)$ --this time there are two possible choices-- 
and consider a family
of rotations~$(U_x, x\in S)$ with the property $U_x e = \n(x)$. 
We can make our construction so
that the maps $x\mapsto \n(x)$ and $x\mapsto U_x$ are continuous a.e.\ in~$S$.

Now consider the random surface crossing times $\htau_n$ of the KSB. We
define this sequence of random variables~$\htau_n$, $n=1,2,3,\ldots$ in the 
following way:
\[
 \htau_1 = \inf\{t\geq 0 : \tX_t\in S\},
\]
and for $n\geq 1$
\[
 \htau_{n+1} = \inf\{t>\htau_n : \tX_t\in S\}.
\]

It is not difficult to prove that for almost all the initial conditions
(as well as e.g.
when the starting location and direction
are chosen independently and uniformly) the surface is intersected only 
finitely many times 
during any finite time interval: $\htau_n\to+\8$ as $n\to\8$\footnote{ 
Indeed, first, analogously to Lemma~\ref{l_abs_continuity}, one can prove
that a.s.\ all the points of crossing are regular points for $S$. 
Then, suppose that for some~$n$
the set $\{\xi_n+s(\xi_{n+1}-\xi_n), s\in[0,1]\}\cap S$ consists of infinitely
many points. In this case, let~$y_0\in S$ be an accumulation point of this set;
it is clear that $\ell_{\xi_n,\xi_{n+1}}\cdot \n(y_0) = 0$. But then, as in
the end of the proof of part~(i) of Lemma~\ref{l_abs_continuity}, 
one can prove that this happens with probability~$0$.}.

Define
\[
 w_n = \left\{
       \begin{array}{ll}
        U_{\tX_{\htau_n}}^{-1}\tV_{\htau_n}, & \text{if } \tV_{\htau_n}\cdot \n(\tX_{\htau_n})\geq 0,\\
        -U_{\tX_{\htau_n}}^{-1}\tV_{\htau_n}, & \text{if } \tV_{\htau_n}\cdot \n(\tX_{\htau_n})<0,\\
       \end{array}
       \right.
\]
i.e., $w_n$ is the relative direction in which the KSB crosses~$S$ for the $n$th time.

By Theorems~\ref{t_conv_KSB} and~\ref{th_gen} we know that~$\tV_t$ has asymptotically uniform
(and independent of the location) distribution, so, if one observes the process when
it is inside some set $A\subset\DD$ with \emph{positive} Lebesgue measure, the
empirical distribution of the direction should be asymptotically uniform as well.
The next result shows that the situation is different when one considers the
instances of intersection with a $(d-1)$-dimensional surface: in this case
the law of intersection is actually the same as the law of reflection
(this result is also connected to Theorem~\ref{t_nonconvex}).

\begin{theo}
\label{t_intersection}
 For a.a.\ initial conditions and any measurable $B\subset\Sph_e$ we have
\begin{equation}
\label{eq_intersection}
\lim_{n\to\8} \frac{1}{n}\sum_{i=1}^n\1{w_i\in B} = \int_B\gamma_d\cos\phi(e,u)\, du.
\end{equation}
\end{theo}

\begin{figure}
\centering
\includegraphics{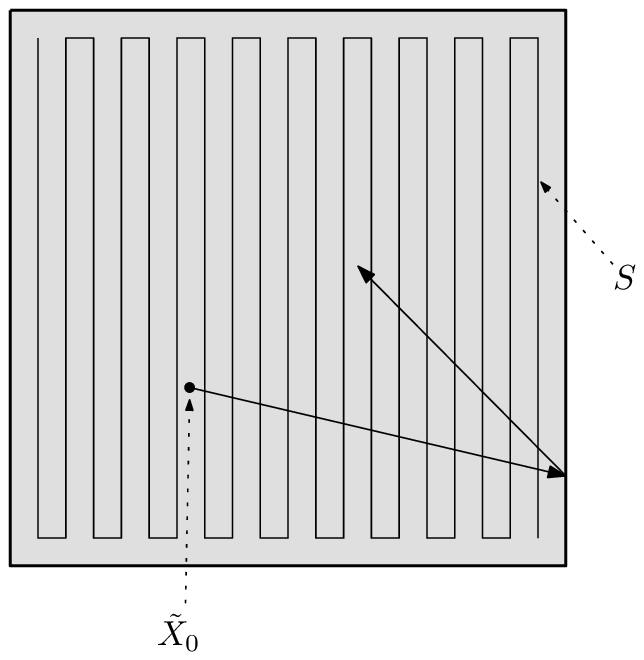}
\caption{With high probability, the (relative) direction of the first crossing of~$S$ is equal
to the initial direction}
\label{fig_cod_barras}
\end{figure}

\noindent
{\bf Remarks.} 
\begin{itemize}
\item From the proof of Theorem~\ref{t_intersection} it can be seen that
a stronger result holds: if we keep track also on the location of crossings, 
then in the long time limit, the 
location is uniformly distributed and is independent 
of the direction.

\item Even if one starts the KSB from the stationary measure, it is generally \emph{not}
true that $(w_n,n\geq 1)$ are identically distributed, and, in particular, it is not true that~$w_1$
has the cosine law. Rather than constructing a counterexample formally, we refer the reader
to Figure~\ref{fig_cod_barras}: it is easy to convince oneself that the law
of~$w_1$ should be close to uniform (and not cosine); more precisely, it converges to
the uniform distribution as the number of vertical components of~$S$ converges
to~$\8$.
\end{itemize}

\section{Some auxiliary results and proof of the 
              existence of the stochastic billiard process}
\label{s_aux}
First, we
recall the definition of the $\alpha$-dimensional Hausdorff measure.
For any Borel set $A\subset\R^d$, define
\[
 \HH^\alpha_\delta (A) = \inf\Big\{\omega_\alpha\sum_{j=1}^\infty
               \Big(\frac{\diam(B_j)}{2}\Big)^\alpha\Big\},
\]
where $\omega_\alpha=\pi^{\alpha/2}/\Gamma(\alpha/2+1)$, and the infimum is
taken over all collections of open sets $(B_j,j\geq 1)$ such that $\diam(B_j)\leq \delta$
for all~$j$ and $A\subset\bigcup_{j=1}^\infty B_j$. Then,
\[
 \HH^\alpha(A) = \lim_{\delta\downarrow 0}\HH^\alpha_\delta (A)
\]
defines the $\alpha$-dimensional Hausdorff measure of the set~$A$.

We need to establish some properties of the domain~$\DD$.

\begin{lmm}
\label{l_basic_Lip}
Suppose that $\diam(\DD)<\infty$. Then we have $|\fr|<\infty$
and the number of connected components of~$\fr$ is finite.
\end{lmm}

\noindent
{\it Proof.} Let us show first that the area of a bounded surface
satisfying the Lipschitz condition is finite. There is a deep
theory related to Lipschitz functions and Hausdorff measure
that can be found in e.g.~\cite{Fed,M}; nevertheless, since we only need
this simple fact, we sketch its proof here.
Consider the following elementary observation: if a function~$f$
satisfies the Lipschitz condition with the constant~$C$ then, for any bounded
$A\subset\Z^{d-1}$,
\begin{eqnarray*}
 \lefteqn{\diam\{z\in\Z^d: (z^{(1)},\ldots,z^{(d-1)})\in A, z^{(d)}=f(z^{(1)},\ldots,z^{(d-1)})\}}
\phantom{****************************}\\
& \leq & (C+1)^{1/2}\diam(A).
\end{eqnarray*}
  From the above definition of the Hausdorff measure it is
elementary to obtain that the area of~$\fr$ is locally finite;
since~$\fr$ is compact, our claim follows.

Let us prove that the number of connected components of~$\fr$ is finite.
Suppose it is not the case; then enumerate these components in some way
and choose any point~$x_j$ from $j$th connected component, $j=1,2,3,\ldots$.
Since~$\DD$ is bounded, there exists $x\in\cDD$ such that~$x$ is an
accumulation point of the sequence $(x_i,i\geq 1)$.
Since~$x_j\in \fr$ and~$\fr$ is a closed set,
we have~$x\in \fr$. But
then, in any arbitrarily small neighborhood of~$x$ one can find pieces of
connected components different from that of~$x$. This contradicts the Lipschitz
condition, since we supposed that in $\BB(x,\eps_x)$ there are no other connected
components of~$\fr$.
\qed

\medskip

For any~$u\in\Sph^{d-1}$, $a,b>0$, define
\[
 \CC_u(a,b) = \Big\{x\in\R^d\setminus\{0\} : \|x\|<a ,
     \tan\phi(u,x/\|x\|)<b\Big\},
\]
so $\CC_u(a,b)$ is a piece of a cone, looking in the direction~$u$.

The next lemma is a technical result that is needed, in particular, to show that
the random variable~$\xi$ is a.s.\ well-defined for all~$n$.
\begin{lmm}
\label{l_abs_continuity}
Suppose that $\diam(\DD)<\infty$.
\begin{itemize}
\item[(i)] Consider any~$x_0\in\cDD$ and let~$B$ be a subset of positive measure of~$\Sph^{d-1}$
such that, for some $\eps_0>0$, we have $x_0+ut\in\DD$ for all $u\in B$, $t\in (0,\eps_0)$.
Let~$\zeta$ be a $B$-valued random variable with the distribution which is
absolutely continuous with respect to the
Haar measure
on~$\Sph^{d-1}$. Then, the distribution of ${\mathsf h}_{x_0}(\zeta)$
is absolutely continuous with
respect to the $(d-1)$-dimensional Hausdorff measure on~$\fr$.
Moreover, with probability~1, 
$\phi_{{\mathsf h}_{x_0}(\zeta)}(\ell_{{\mathsf h}_{x_0}(\zeta),x_0})<\pi/2$.
\item[(ii)] 
The conditional distribution  of~$\xi_{n+1}$ given $\xi_n=x_1\in\RR$,
is absolutely continuous with
respect to the $(d-1)$-dimensional Hausdorff measure on~$\fr$.
\end{itemize}
\end{lmm}

\noindent
{\it Proof.}
First, we prove the part~(i). We have to prove that the measure 
${\tilde \mu}(A)=\PP[{\mathsf h}_{x_0}(\zeta)\in A]$ is absolutely continuous with
respect to the $(d-1)$-dimensional Hausdorff measure on~$\fr$. 
Consider any~$A\subset\fr$ such that $|A|=0$.
This means that, for any~$\eps>0$, there are sequences $(z_i\in\R^d,i\geq 1)$ and
$(r_i>0,i\geq 1)$ such that 
\[
 A\subset\bigcup_{i=1}^\infty\BB(z_i,r_i) \quad \text{ and }\quad \sum_{i=1}^\infty r_i^{d-1}<\eps.
\]
Without restricting the generality, one can suppose that the distance from~$x_0$ to any
of $\BB(z_i,r_i)$ is at least $\eps_0/2$.
Let
\[
 B_i = \{u\in\Sph^{d-1} : \text{ there exists $t>0$ such that }x_0+ut\in\BB(z_i,r_i)\}
\]
be the projection of $\BB(z_i,r_i)-x_0$ on~$\Sph^{d-1}$.
Clearly, there exists $C_{\eps_0}>0$ such that $|B_i|\leq C_{\eps_0} r_i^{d-1}$.
If $g(u)$ is the density of~$\zeta$, we can write 
\[
 \PP[{\mathsf h}_{x_0}(\zeta)\in A] \leq \int_{\cup_{i=1}^\infty B_i} g(u)\, du;
\]
on the other hand, 
\[
 \Big|\bigcup_{i=1}^\infty B_i\Big| \leq \sum_{i=1}^\infty |B_i| 
  \leq C_{\eps_0}\sum_{i=1}^\infty r_i^{d-1} < \eps C_{\eps_0}.
\]
Since~$\eps$ is arbitrary, by the continuity property of the Lebesgue integral
we obtain that ${\tilde \mu}(A)=0$. 

To complete the proof of part~(i), observe 
that, with  ${\hat \gamma}$ the density of~$\zeta$,
\begin{eqnarray*}
 \PP[{\mathsf h}_{x_0}(\zeta)\in A] 
&=&
\int_{\Sph^{d-1}}  {\hat \gamma}(v) \1{{\mathsf h}_{x_0}(v)\in A} \, dv\\
&=&
\int_{\Sph^{d-1}}  {\hat \gamma}(v) \1{{\mathsf h}_{x_0}(v)\in A \cap \RR} 
\, dv\\
&=& 
 \int_A {\hat \gamma}(\ell_{x_0,y})
             \frac{\cos\phi_y(\ell_{y,x_0})}{\|x_0-y\|^{d-1}}
                \1{y\in\RR,x_0\leftrightarrow y}\, dy.
\end{eqnarray*}
The second line is from the above absolute continuity,
and the last line is obtained by the
change of variables from
 $v \in  \Sph^{d-1}$ to $y={\mathsf h}_{x_0}(v) \in \fr$,
which can be performed at regular points of the boundary provided
that $\phi_y(\ell_{y,x_0})<\pi/2$.

Finally, it is an easy exercise to deduce the part~(ii) from the part~(i)
(decompose $\Sph_{\n(x_1)}$ into the union of $D_j:=\{u\in\Sph_{\n(x_1)} :
{\mathsf r}_{x_1}(u)\in[1/(j+1),1/j)\}$, $j\geq 1$, and $D_0:=\{u\in\Sph_{\n(x_1)} :
{\mathsf r}_{x_1}(u)\geq 1\}$, and then use the part~(i)).
%
\qed

\medskip

Define, for all~$x\in\fr$,
\begin{equation}
\label{def_alpha_x}
 {\tilde\alpha}_x = \sup\{s>0 : \text{ there exists }u\in\Sph^{d-1} 
     \text{ such that } x + \CC_u(s,s) \subset \DD\}.
\end{equation}
Let us prove that, under the Lipschitz condition
formulated in Section~\ref{s_gen_notation}, the domain~$\DD$ has the following properties:

\begin{lmm}
\label{good_domain}
\begin{itemize}
\item[(i)] For any compact $A\subset \R^d$ there exists a constant~$q(A)>0$ such that
${\tilde\alpha}_x \geq q(A)$ for all $x\in \fr \cap A$.
\item[(ii)] Suppose that $\diam(\DD)<\infty$. Then,
for any $x\in\fr$, there exist $\eps,\eps',\delta>0$ and $y_x\in\RR$ such that
$K(y,z)>\delta$ for all $z\in\BB(x,\eps)\cap\RR$, $y\in\BB(y_x,\eps')$.
\item[(iii)] Suppose that $x,y\in\RR$ are such that $K(x,y)>0$. Then there 
exist $\eps,\eps',\delta>0$ such that $K(x',y')>\delta$
for all $x'\in\BB(x,\eps)$, $y'\in\BB(y,\eps')$.
\end{itemize}
\end{lmm}

\noindent
{\it Proof.}
Consider an arbitrary $x\in\fr$ and fix~$a_0<1$ in such a way that
 $a_0(1-a_0^2)^{-1/2}<(2L_x)^{-1}$ (recall that we assumed $L_x>1$).
 From the Lipschitz condition we obtain that 
for any $y\in\fr\cap\BB(x,a_0\eps_x)$ it holds that
$y+\CC_{{\mathfrak I}_x^{-1}(e_d)}((1-a_0)\eps_x,L_x^{-1})\subset(\DD\cap \BB(x,\eps_x))$.
The collection $(\fr\cap\BB(x,a_0\eps_x), x\in \fr \cap A)$
is a covering of $\fr \cap A$; choose a finite subcovering to establish 
the property~(i).

\begin{figure}
\centering
\includegraphics{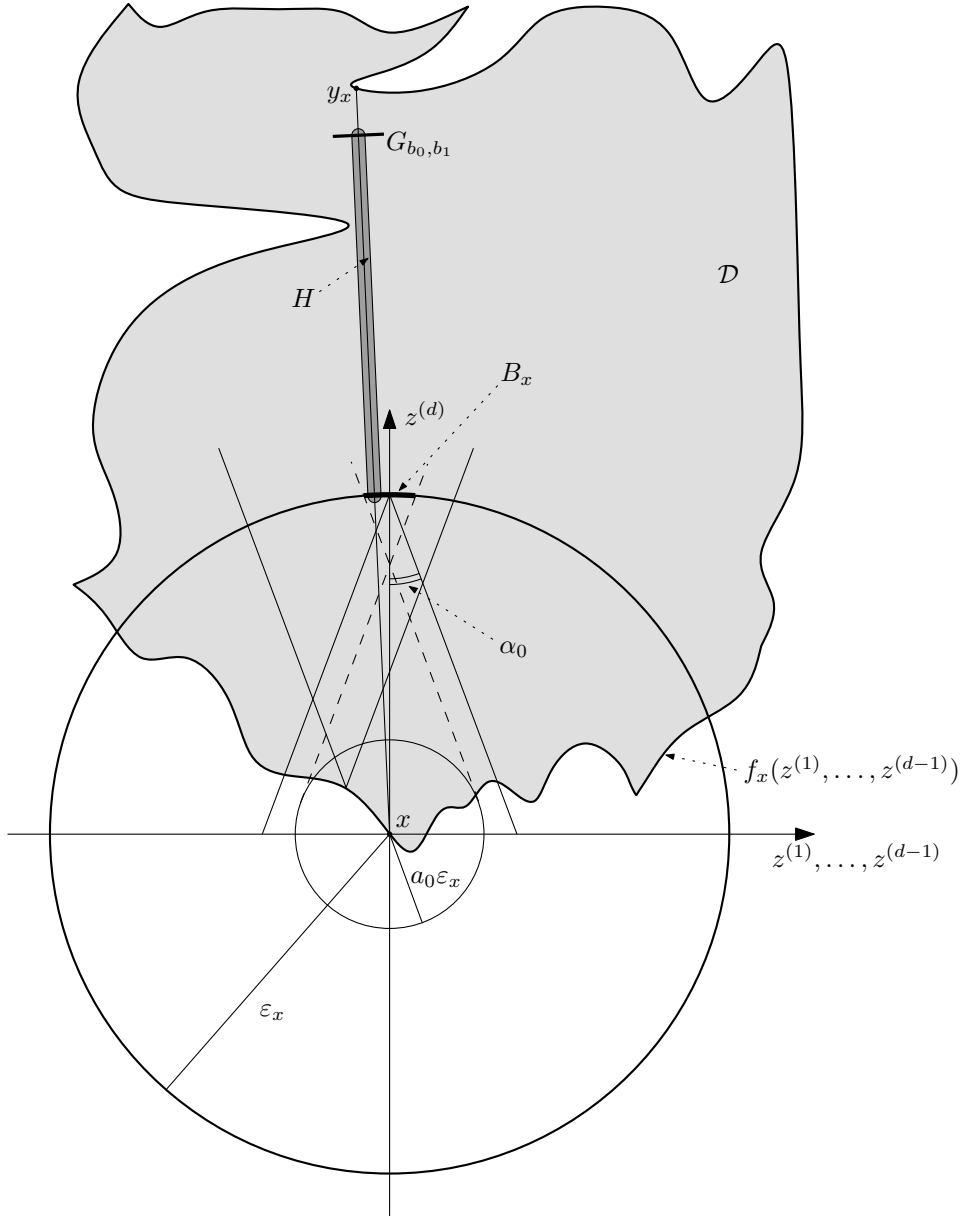}
\caption{On the proof of Lemma~\ref{good_domain} (ii), $\alpha_0:=\arctan(2L_x)^{-1}$}
\label{fig_cones}
\end{figure}

Let us prove the property~(ii). Suppose that $\diam(\DD)<\infty$ and observe that
one can obtain (see Figure~\ref{fig_cones})
\[
 \DD\cap \BB(x,\eps_x) \cap
\bigcap_{y\in\fr\cap\BB(x,a_0\eps_x)}(y+\CC_{{\mathfrak I}_x^{-1}(e_d)}(\infty,(2L_x)^{-1}))
 \neq \emptyset.
\]
This set being open, this implies that $|B_x|>0$ (here $|\cdot|$ stands 
for the $(d-1)$-dimensional Hausdorff measure), where
\[
 B_x = \{z:\|z-x\|=\eps_x\}\cap\bigcap_{y\in\fr\cap\BB(x,a_0\eps_x)}
     (y+\CC_{{\mathfrak I}_x^{-1}(e_d)}(\infty,(2L_x)^{-1})).
\]
Note that any point of~$B_x$ is seen from any point of $\fr\cap\BB(x,a_0\eps_x)$, 
and for all $y\in B_x$ it holds that
\begin{equation}
\label{ugol_xy}
 \phi_x(\ell_{x,y}) \geq \arctan L_x^{-1} - \arctan(2L_x)^{-1}.
\end{equation}
Moreover, from Lemma~\ref{l_abs_continuity}~(i) it follows that
 one can find a point~$y_x\in\RR$ such that $\eps_x\ell_{x,y_x}\in B_x$,
$\phi_{y_x}(\ell_{y_x,x})<\pi/2$, and $x\leftrightarrow y_x$.

For $b_1,b_2>0$ denote 
\[
 G_{b_1,b_2} = \{z\in\R^d: (z-y_x-b_1\ell_{y_x,x})\cdot \ell_{y_x,x} = 0, 
            \|z-y_x-b_1\ell_{y_x,x}\|<b_0\}.
\] 
Since~$y_x\in\RR$ and $\phi_{y_x}(\ell_{y_x,x})<\pi/2$,
one can choose positive and small enough $b_1,b_2,\eps'$ (also with the
property $b_1<\|x-y_x\|-\eps_x$) such that
\begin{itemize}
\item $G_{b_1,b_2} \subset \DD$;
\item $\fr\cap\BB(y_x,\eps') \subset \RR$;
\item $y\leftrightarrow z$ for all $z\in G_{b_1,b_2}$, $y\in\fr\cap\BB(y_x,\eps')$.
\end{itemize}

Clearly, there is a small enough $b_2\leq b_0$ such that for the set
\[
 H = \{z: \text{ there exists } t\in[b_1,\|x-y_x\|-\eps_x] 
            \text{ such that } \|z-t\ell_{y_x,x}\|\leq b_2\}
\]
it holds that $H\subset\DD$ and $\{z\in H : \|x-z\|=\eps_x\}\subset B_x$.
Then, one can find positive and small enough $s,\eps'_1\leq \eps', a_1\leq a_0$
such that for any $y\in\fr\cap\BB(y_x,\eps'_1)$
\begin{itemize}
\item $(y+\CC_{\ell_{y_x,x}}(\infty,s))\cap \{z:\|z-x\|=\eps_x\} \subset B_x$;
\item $(y+\CC_{\ell_{y_x,x}}(\infty,s))\cap 
      \{z: (z-y_x-b_1\ell_{y_x,x})\cdot\ell_{y_x,x}\}\subset G_{b_0,b_1}$;
\item $\BB(x,a_1\eps_x) \subset (y+\CC_{\ell_{y_x,x}}(\infty,s))$.
\end{itemize}
Then, by construction, we obtain that $x'\leftrightarrow y$ for any
$x'\in\fr\cap\BB(x,a_1\eps_x)$, $y\in\fr\cap\BB(y_x,\eps'_1)$.
For part~(ii) it remains only to prove that the above parameters can be chosen
in such a way that $K(x',y)$ is uniformly positive for regular points
from the neighborhoods of~$x$ and~$y_x$. This can be obtained easily,
using that $\phi_{y_x}(\ell_{y_x,x})<\pi/2$ and~(\ref{ugol_xy}).

Proving part~(iii) is very similar to the last part of the proof of part~(ii),
so we do not provide the details.
\qed

\medskip

Now, we are ready to prove that the stochastic billiard process~$X_t$
is well defined for all~$t\geq 0$.

\medskip
\noindent
{\it Proof of Proposition~\ref{tau_to_infty}.}
Fix any $T>0$ and let us prove that 
\begin{equation}
\label{lim<T}
\lim_{n\to\infty}\tau_n > T \qquad \text{a.s.}
\end{equation}
Consider the compact set $A=\{x\in\R^d : \|x-x_0\|\leq T\}$ 
(recall that $\xi_0=x_0$); by Lemma~\ref{good_domain}~(i), 
we have that ${\tilde \alpha}_x \geq q(A)>0$ 
for all $x\in \fr \cap A$, where ${\tilde \alpha}_x$
was defined in~(\ref{def_alpha_x}). 
Denote $B_u=\{v\in\Sph^{d-1} : \tan\phi(u,v)<q(A)\}$.
Then, from~(\ref{support_halfsph}) it is elementary to obtain that 
\begin{equation}
\label{prob_ubezhat'}
  \PP[\|\xi_{n+1} - \xi_n\| \geq {\tilde\alpha}_x \mid \xi_n=x]
                \geq \inf_{u: B_u\subset\Sph_e}\PP[\eta\in B_u] > 0
\end{equation}
for any $x\in \fr \cap A$.
This means that, as long as $\tau_n\leq T$ (and, consequently, $\xi_n\in \fr\cap A$),
there is a uniformly positive probability 
that $\tau_{n+1}-\tau_n = \|\xi_{n+1} - \xi_n\| \geq q(A)$.
This proves~(\ref{lim<T}), and, since~$T$ is arbitrary, Proposition~\ref{tau_to_infty}
holds.
\qed

\section{Proofs for the random walk}
\label{s_proofs_KRW}
Here we prove the results concerning the discrete-time random walk~$\xi_n$.

\medskip
\noindent
{\it Proof of Theorem~\ref{t_conv_KRW}.}
First, our goal is to prove that the D\"oblin condition holds: 
there exist~$n_0$, ${\hat \eps}>0$ such
that for all $x,y\in\RR$ 
\begin{equation}
\label{Doeblin}
 K^{n_0}(x,y) \geq {\hat \eps}.
\end{equation}

To prove~(\ref{Doeblin}), let us first suppose that~$x$ and~$y$ are in the same
connected component of~$\fr$. 
Now, from Lemma~\ref{good_domain} it follows that
for any $x\in\fr$ we can find $\eps(x),\delta(x)>0$, $s(x)\in\RR$ such that
\begin{equation}
\label{horosho_vidno}
\inf_{\substack{y\in\RR\cap\BB(x,\eps(x)),\\ z\in\RR\cap\BB(s(x),\delta(x))}} K(y,z) > 0.
\end{equation}

Note that $(\fr\cap\BB(z,\eps(z)), z\in\fr)$ is a covering of $\fr$, so let us choose
a finite subcovering $(\fr\cap\BB(x_i,\eps(x_i)), i=1,\ldots,N)$. 
 Abbreviate $B_{i,j}:=\fr\cap\BB(x_i,\eps(x_i))\cap\BB(x_j,\eps(x_j))$, let
\begin{eqnarray*}
 \beta &=& \min_{i=1,\ldots,N} 
    \inf_{\substack{y\in\RR\cap\BB(x_i,\eps(x_i)),\\ z\in\RR\cap\BB(s(x_i),\delta(x_i))}} K(y,z),\\
 \theta_1 &=& \min_{i,j: B_{i,j}\neq\emptyset} |B_{i,j}|,\\
 \theta_2 &=& \min_{i=1,\ldots,N} |\fr\cap\BB(s(x_i),\delta(x_i))|
\end{eqnarray*}
(note that all the quantities, defined above, are strictly positive).

Since $x,y$ are from the same connected component, there exists $m\leq N$ 
and $k(1),\ldots,k(m)\in\{1,\ldots,N\}$, 
such that $x\in \RR\cap\BB(x_{k(1)})$, $y\in \RR\cap\BB(x_{k(m)})$,
and $B_{k(i),k(i+1)}\neq \emptyset$ for all $i=1,\ldots,m-1$.
Note also that one can suppose that~$m$ is the same for all~$x,y$
(and even that it does not depend on the connected component
of~$\fr$ to which belong~$x$ and~$y$), since we
formally can choose $k(i)=k(i+1)$ in the above sequence. 
So, abbreviating $D_j=\RR\cap\BB(s(x_{k(j)}),\delta(x_{k(j)}))$, we can write
\begin{eqnarray*}
 K^{2m}(x,y) &\geq & \int_{D_1} dz_1 \int_{B_{k(1),k(2)}}dw_1 
      \int_{D_2}dz_2\\
     && \ldots \int_{D_{m-1}}dz_{m-1} 
      \int_{B_{k(m-1),k(m)}}dw_{m-1}
        \int_{D_m}dz_m\\
     &&  K(x,z_1)K(z_1,w_1) K(w_1,z_2)\ldots K(w_{m-2},z_{m-1})\\
     &&{}\times K(z_{m-1},w_{m-1}) K(w_{m-1},z_m)K(z_m,y)\\
   &\geq &\beta^{2m-1}\theta_1^{m-1}\theta_2^m,
\end{eqnarray*}
which proves~(\ref{Doeblin}) for the case when $x,y$ are in the same connected
component of~$\fr$.

\begin{figure}
\centering
\includegraphics{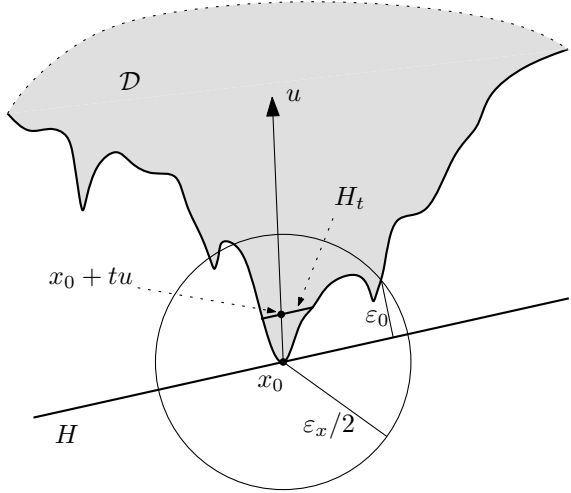}
\caption{Finding two regular points which see each other}
\label{fig_exposed}
\end{figure}

Now, we sketch the proof of~(\ref{Doeblin}) for the case when $x,y$ are in different connected
components of~$\fr$; we leave the details to the reader. Suppose that there are~$M$
connected components of~$\fr$, so that $\fr=\cup_{i=1}^M G_i$. Let $a_{ij}=1$
if there exist $z_{ij} \in G_i\cap\RR$, $z_{ji}\in G_j\cap\RR$ 
such that $K(z_{ij},z_{ji})>0$, and $a_{ij}=0$ otherwise. In words,
$A=(a_{ij})_{i,j=1\ldots,M}$ is the matrix that describes the
``graph'' of the set of connected components of~$\fr$.
It is clear that the matrix~$A$ is irreducible,
let us prove that it is aperiodic. 
Let $G_{i_0}$ be the 
connected component of~$\fr$ which is the boundary of the unique infinite
connected component of $\R^d\setminus \DD$, we show that $a_{i_0i_0}=1$.
Let~$x_0$ be an exposed point of the convex hull of~$\cDD$\footnote{i.e.,
such that there exists an hyperplane $H$ which intersects the convex 
hull of~$\cDD$ only at $\{x_0\}$}, and let~$H$ be
a supporting hyperplane of the convex hull of~$\cDD$ at~$x_0$, i.e., 
$H\cap\cDD = \{x_0\}$. From the Lipschitz condition (see also the proof
of Lemma~\ref{good_domain}~(i)) we know that for some $u\in\Sph^{d-1}, b_1,b_2>0$,
it holds that $\CC_u(b_1,b_2)\subset\DD$. Let (see Figure~\ref{fig_exposed})
\[
 H_t = \{tu+z : z\in H \text{ is such that } tu+\alpha z \in\DD \text{ for all }\alpha\in[0,1)\}.
\]
Let
\[
 \eps_0 = \inf\{\|x-y\| : x\in H, y\in\fr
, \|y-x_0\|=\eps_x/2
\} > 0.
\]
Clearly, $H_t$ is a bounded
set for all $t\in (0,\eps_0)$ and $\overline{H_t}\cap\fr$ contains only the
points of the connected component of~$x_0$.
Let~$\zeta_1$ be a random variable with uniform distribution on~$(0,\eps_0)$,
and~$\zeta_2$ a random variable with uniform distribution on~$\Sph^{d-1}\cap (H-x_0)$.
Define $y_1\in H_{\zeta_1}$ such that 
$y_1=\zeta_1 u + \alpha_1\zeta_2$ for some $\alpha_1>0$,
and $y_2\in H_{\zeta_1}$ such that
 $y_2=\zeta_1 u - \alpha_2\zeta_2$ for some $\alpha_2>0$.
By construction, we have that $y_1,y_2$ are in the connected
component of~$x_0$ and $y_1\leftrightarrow y_2$.
Similarly to the proof of Lemma~\ref{l_abs_continuity}~(i), one can show that
a.s.\ $y_1,y_2\in\RR$ and $\phi(\n(y_1),\ell_{y_1,y_2})\in[0,\pi/2)$,
$\phi(\n(y_2),\ell_{y_2,y_1})\in[0,\pi/2)$.

Thus, there exists~$m_1$ such that $A^{m_1}>0$.
Consider any $x\in G_i$, $y\in G_j$, there exist $i=k_0,k_1,\ldots,k_{m_1-1},k_{m_1}=j$
such that $a_{k_lk_{l+1}}=1$ for $l=0,\ldots,m_1-1$. Then, roughly speaking,
we go from~$x$ to $z_{k_0k_1}$ in~$2m$ steps, then 
(recall Lemma~\ref{good_domain}~(iii)) jump from $z_{k_0k_1}$
to $z_{k_1k_0}$, then go from $z_{k_1k_0}$ to $z_{k_1k_2}$ in~$2m$ steps
etc. This shows that~(\ref{Doeblin}) holds with $n_0=2m(m_1+1)$.

Now, having~(\ref{Doeblin}) at hand, it is straightforward to observe that
there exists a coupling of two versions~$\xi$ and~$\xi'$ of the random walk
with arbitrary starting conditions such that $\PP[\xi_{n_0}=\xi'_{n_0}]\geq 
{\bar \eps}:={\hat\eps}|\fr|$.
This in turn means that one can construct the two processes~$\xi$ and~$\xi'$
together with a random variable~${\tilde T}$ in such a way that 
$\xi_{{\tilde T}+k} = \xi'_{{\tilde T}+k}$ for all~$k\geq 0$, and this random
variable~${\tilde T}$ (which is the \emph{exact coupling time} for the two processes)
has the property
\[
 \PP[{\tilde T} > k n_0] \leq (1-{\bar \eps})^k.
\]
Applying Theorem~5.1 of Chapter~4 from~\cite{T}, we obtain~(\ref{eq_conv_KRW}).

The absolute continuity of~${\hat\mu}$ with respect to the $(d-1)$-dimensional
Hausdorff measure follows from Lemma~\ref{l_abs_continuity}~(ii), and 
the equation~(\ref{eq_stat_dens}) follows from the fact that if the distribution of~$\xi_0$
is~${\hat\mu}$, then so is the distribution of~$\xi_1$. Note that
$\psi(x)=\int_{\fr}\psi(y)K^{n_0}(y,x)\, dy$ for $x \in \RR$, 
so from~(\ref{Doeblin}) it
follows that $\inf_{\fr}\psi \geq {\hat\eps}$.

\medskip

We end by proving~(\ref{eq_CLT_KRW}). This is now completely standard, 
a direct consequence of 
(ii)-(iv), Theorem 17.0.1 in~\cite{MT}.
\qed

\section{Proofs for the continuous time stochastic billiard}
\label{s_proofs_KSB}
In this section we prove the results concerning the continuous-time 
stochastic billiard.

\medskip

\noindent
{\it Proof of Proposition~\ref{prop:pbmg}.}
For $0\leq u \leq t$, we have
\begin{eqnarray} 
\lefteqn{\E\Big[f(X_t,V_t)-f(X_u,V_u)- \int_u^t  \LL f(X_s,V_s)ds \,\Big|\, \FF_u\Big]}\\
\nonumber
&=& \E\Big[ \sum_{n \geq 0} {\bf 1}\{u<\tau_n,\tau_{n-1}  \leq t\} 
\Big(f(X_{\tau_n \wedge t},V_{\tau_n \wedge t})-
f(X_{\tau_{n-1} \vee u},V_{\tau_{n-1} \vee u}) \\
&& \qquad \qquad   \nonumber
- \int_{\tau_{n-1} \vee u}^{\tau_n \wedge t}  
\LL f(X_s,V_s)ds\Big) \,\Big|\, \FF_u\Big]\\
&=& \sum_{n \geq 0} \E\Big[ {\bf 1}\{u<\tau_n, \tau_{n-1}  \leq t\} 
\Big( f(X_{\tau_n \wedge t},V_{\tau_n \wedge t})-
f(X_{\tau_{n-1} \vee u},V_{\tau_{n-1} \vee u})  \nonumber\\
&& \qquad \qquad -\int_{\tau_{n-1} \vee u}^{\tau_n \wedge t}  
\LL f(X_s,V_s)ds\Big) \,\Big|\, \FF_u\Big] \label{eq:f1}\\
&=& \sum_{n \geq 0} \E\Big[ {\bf 1}\{u<\tau_n \leq t\} 
\Big(f(X_{\tau_n},V_{\tau_n})-
f(X_{\tau_{n}},V_{\tau_{n-1}}) 
\Big) \,\Big|\, \FF_u\Big]   \label{eq:f2}\\
&=&0, \label{eq:f3}
\end{eqnarray}
where we have used Fubini's theorem in (\ref{eq:f1}), the fundamental theorem 
of calculus in (\ref{eq:f2}), and the boundary conditions (\ref{eq:dom}) in
(\ref{eq:f3}). This proves that~(\ref{eq:pbmart}) defines a martingale.
\qed

\medskip

\medskip
\noindent
{\it Proof of Theorem~\ref{t_conv_KSB}.}
To prove~(\ref{eq_conv_KSB}), we use the approach of~\cite{T}.
By compactness of the state space $\DD \times \Sph^{d-1}$, there exists 
at least one invariant probability measure $\chi$ for the Markov process
 $(X_t,V_t)$.
Consider the process $(X_t,V_t)$ starting from an arbitrary initial 
position~$x_0$
and initial direction~$v_0$, and consider also the stationary version 
$({\hat X}_t,{\hat V}_t)$ of this process (i.e., 
$({\hat X}_0,{\hat V}_0) \sim \chi$, which implies that
$({\hat X}_t,{\hat V}_t) \sim \chi$ for all $t$).
Our goal now is to construct an exact coupling of
$(X_t,V_t)$ and $({\hat X}_t,{\hat V}_t)$, i.e., we construct these processes
together with a random variable~${\hat T}$ (the coupling time) on a same probability space
in such a way that 
\[
 (X_{{\hat T}+s},V_{{\hat T}+s}) = ({\hat X}_{{\hat T}+s},{\hat V}_{{\hat T}+s})
\]
for all~$s\geq 0$. Also, we show that~${\hat T}$ has exponential tail, which in its turn
ensures~(\ref{eq_conv_KSB}).

\begin{df}
\label{def_alpha_cont}
We say that a pair of random variables $(X,T)$, $X\in\fr,T\in\R_+$ is $\alpha$-continuous
on the set $A\times B$, $A\subset\fr, B\subset\R_+$, if for any 
measurable $A_1\subset A$, $B_1\subset B$
\[
 \PP[X\in A_1, T\in B_1] \geq \alpha |A_1||B_1|.
\]
\end{df}

For any $x\in\fr$, $\eps>0$, define $\Delta_x^\eps=\BB(x,\eps)\cap\RR$. 
In the sequel, we will write $\xi^a_n$ for the position of the random walk that 
started in~$a$, $\tau^a_n$ for the corresponding local times.
Note that, if we consider the pair $(\xi^a_1,\tau^a_1)$,
its joint distribution is singular with respect
to the Lebesgue measure (although typically the distribution of each component is absolutely
continuous). The next lemma shows that in two steps the situation is already better.

Fix an arbitrary $x_1\in\RR$ with the following property: there exists $y_1\in\RR$
such that $K(x_1,y_1)>0$ and $\n(y_1)\cdot \ell_{y_1,x_1} < 1$.
\begin{lmm}
\label{small_dens}
There exist $\alpha_1,\eps>0$, ${\hat r}_2>{\hat r}_1>0$ such that,
for any $x\in \Delta_{x_1}^\eps$, the pair $(\xi^x_2,\tau^x_2)$
is $\alpha_1$-continuous on $\Delta_{x_1}^\eps\times ({\hat r}_1,{\hat r}_2)$.
\end{lmm}

\noindent
{\it Proof.}
Choose small enough~$\eps,\eps'$ in such a way that
$K(x,y)>h_1$ for all 
$x\in\Delta_{x_1}^\eps$, $y\in\Delta_{y_1}^{\eps'}$ (this is possible 
by Lemma~\ref{l_abs_continuity}~(iii); note that
$\Delta_{x_1}^\eps\cup \Delta_{y_1}^{\eps'}\subset\RR$).
Denote by 
\[
H=\{y'\in\R^d:\; (y'-y_1)\cdot\n(y_1)=0\}
\]
the tangent hyperplane at the point~$y_1$, and let
\[
v_1=(\n(y_1)\cdot \ell_{y_1,x_1})\n(y_1)-\ell_{y_1,x_1}.
\]
Note that
\begin{itemize}
\item $v_1+y_1\in H$ (in fact, $v_1$ is the projection of the vector $\n(y_1)-\ell_{y_1,x_1}$
onto the hyperplane $H-y_1$);
\item $v_1\cdot \ell_{y_1,x_1}<0$ (since 
 $v_1\cdot \ell_{y_1,x_1}=(\n(y_1)\cdot \ell_{y_1,x_1})^2-1$ is negative by the choice of $y_1$).
\end{itemize}

For notational convenience, for the rest of the proof we assume that $d\geq 3$;
it is straightforward to adapt the proof for the case $d=2$ (in fact, the proof in
this case is much simpler, since the sets $G,G_1,G_1^+,G_1^-$ defined below
would be simply one point sets).

\begin{figure}
\centering
\includegraphics{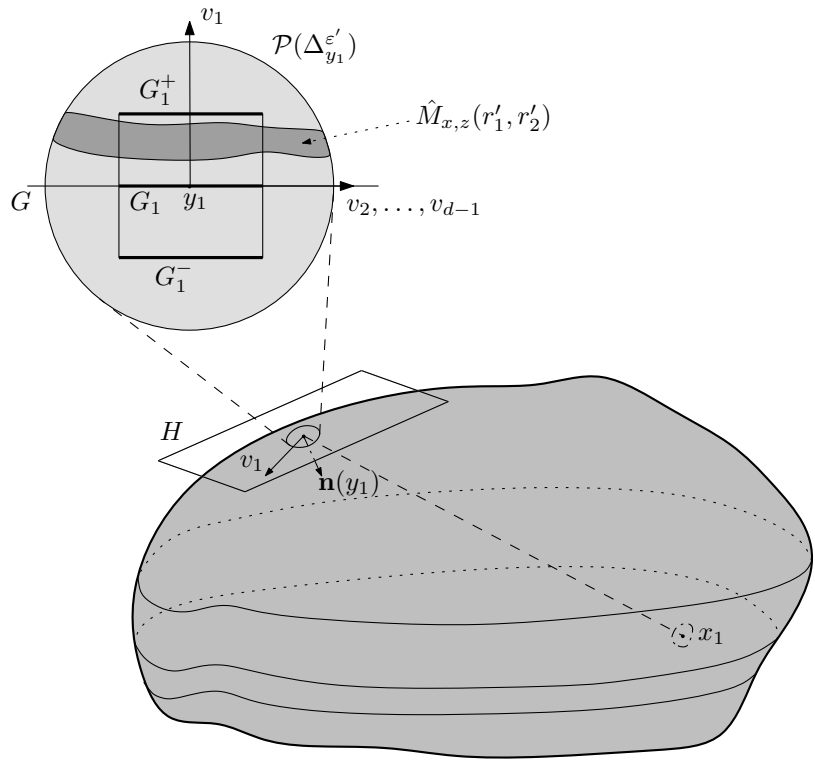}
\caption{On the proof of Lemma~\ref{small_dens}}
\label{fig_alpha_cont}
\end{figure}

Consider the $(d-2)$-dimensional plane (see Figure~\ref{fig_alpha_cont})
\[
 G=\{z\in\R^d:\; (z-y_1)\cdot v_1=0\}\cap H.
\]
One can easily see that if $y'-y_1\in G$, then $(y'-y_1)\cdot \ell_{y_1,x_1}=0$.
Let us choose an orthonormal basis $v_2, \ldots, v_{d-1}$ in~$G-y_1$, so that
$v_1, v_2, \ldots, v_{d-1}$ is an orthonormal basis in $H-y_1$. 
Let ${\mathcal P}: \partial D\to H$ be the projection operator
defined by
\[
 {\mathcal P}y = y-\big((y-y_1)\cdot\n(y_1)\big)\n(y_1),
\]
and denote, with $R^x(y,z) := \|x-y\|+\|y-z\|$,
\[
\theta_v(t):=R^{x_1}({\mathcal P}^{-1}(y_1+vt),x_1)=2\|{\mathcal P}^{-1}(y_1+vt)-x_1\|.
\]
As ${\mathcal P}^{-1}(y_1+vt)=y_1+vt+O(t^2)$ when $t\to 0$, using 
also the fact that $\nabla\|x\|=\frac{x}{\|x\|}$, we have
\begin{equation}
\label{der_vector}
\theta'_v(0)=\frac{2(y_1+vt-x_1)}{\|y_1+vt-x_1\|}\cdot v\Big|_{t=0}=2\ell_{x_1,y_1}\cdot v
\end{equation}
(note that $y_1-x_1=-\ell_{y_1,x_1}\|y_1-x_1\|$).

Since $v_1, \ldots, v_{d-1}$ is a basis in $H-y_1$, any
$s\in H-y_1$ can be uniquely represented as
\[
s=s^{(1)}v_1+\cdots +s^{(d-1)}v_{d-1},
\]
so we regard $s^{(1)},\ldots,s^{(d-1)}$ as the coordinates of~$s$.
 Note that ${\mathcal P}^{-1}(s+y_1)\in \fr$
and let us consider the function
\[
\Psi_{x,z}(s)=R^x({\mathcal P}^{-1}(s+y_1), z).
\]
 From~(\ref{der_vector}) we obtain (recall that~$G$ is perpendicular to $\ell_{x_1,y_1}$)
\begin{eqnarray}
\label{psi1}
\frac{\partial \Psi_{x_1, x_1}(s)}{\partial s^{(1)}}
   \Big|_{s=0}&=&2\ell_{x_1,y_1}\cdot v_1=:h>0,\\
\frac{\partial \Psi_{x_1, x_1}(s)}{\partial s^{(i)}}\Big|_{s=0}&=&0, 
                 \quad\text{ for }i=2, \ldots, d-1.\nonumber
\end{eqnarray}
Since $\Psi_{x,z}(s)$ is continuously differentiable, for small enough $\eps,\eps'$
we obtain that for $x, z\in \Delta_{x_1}^\eps$ 
and $s+y_1\in \Delta_{y_1}^{\eps'}$,
\begin{eqnarray}
\label{psi2}
\frac{\partial \Psi_{x, z}(s)}{\partial s^{(1)}}&\in & \Big[\frac{3h}{4},2h\Big],\\
\label{psi3}
\frac{\partial \Psi_{x_1, x_1}(s)}{\partial s^{(i)}}&\leq &\frac{h}{4}, 
                 \quad\text{ for }i=2, \ldots, d-1.
\end{eqnarray}

Denote $G_1=\{z\in G: \|z-y_1\|\le \eps'/2\}$, 
$G_1^+=G_1+\frac{\eps'}{2}v_1$, $G_1^-=G_1-\frac{\eps'}{2}v_1$
(again, see Figure~\ref{fig_alpha_cont}). Clearly, if~$\eps'$
is small enough,
\[
 \{s+tv_1: s\in G_1, t\in[-\eps'/2,\eps'/2]\}+y_1 \subset {\mathcal P}(\Delta_{y_1}^{\eps'}).
\]
Let 
\[
\hat r_1=\sup_{x,z\in\Delta_{x_1}^\eps}\sup_{s\in G_1^-}\Psi_{x,z}(s),\quad
\hat r_2=\inf_{x,z\in\Delta_{x_1}^\eps}\inf_{s\in G_1^+}\Psi_{x,z}(s).
\]
Let us show that $\hat r_1<\hat r_2$. Indeed,
\[
\Psi_{x,z}\Big(\frac{\eps'}{2}v_1\Big)-\Psi_{x,z}\Big(-\frac{\eps'}{2}v_1\Big)
         \geq \frac{3h}{4}\eps'
\] 
by~(\ref{psi2}), and for any $s\in G_1^+$
\[
\Psi_{x,z}\Big(\frac{\eps'}{2}v_1\Big)-\Psi_{x,z}(s) \leq \frac{h}{4}\eps'
\]
by (\ref{psi3}), and the similar inequality holds for any $s\in G_1^-$.
So, we obtain that 
${\hat r}_2-{\hat r}_1\geq \frac{3h}{4}\eps' - \frac{h}{4}\eps' - \frac{h}{4}\eps' = \frac{h}{4}\eps'$.

Let
\begin{eqnarray*}
M_{x,z}(r'_1,r'_2) &=& \{y\in\Delta_{y_1}^{\eps'} : R^x(y,z)\in (r'_1,r'_2)\},\\
{\hat M}_{x,z}(r_1', r_2') &=&
       \{s\in {\mathcal P}(\Delta_{y_1}^{\eps'}) - y_1:\; \Psi_{x,z}(s)\in(r_1', r_2')\};
\end{eqnarray*}
observe that $M_{x,z}(r'_1,r'_2)={\mathcal P}^{-1}({\hat M}_{x,z}(r_1', r_2'))$.

Note the following simple fact: if for all $x\in[a_1, a_2]$ we have $c_1<g'(x)<c_2$
where $0<c_1<c_2<\infty$, then for any $I\subset[f(a_1), f(a_2)]$ it holds that
$c_2^{-1}|I|\le |f^{-1}(I)|\le c_1^{-1}|I|$.
Using this for $g_s(t)=\Psi_{x,z}(s+v_1 t)$ (recall that, by~(\ref{psi2}), 
$g'(t)\in[3h/4, 2h]$) we obtain 
\[
|\hat M_{x,z}(r_1', r_2')| \geq \int_{G_1}\big|g_s^{-1}\big((r'_1,r'_2)\big)\big|\,ds
    \geq |G_1| (2h)^{-1} (r_1'- r_2')
\]
for all
$(r_1', r_2')\subset ({\hat r}_1,{\hat r}_2)$. Clearly, for small enough $\eps'$, there exists 
a constant~$c'$ such that for any measurable $A\subset {\mathcal P}^{-1}(\Delta_{y_1}^{\eps'})$
we have $|{\mathcal P}^{-1}(A)|\ge c'|A|$.
This implies that
there exists $\alpha'>0$ such that for all $x,z\in\Delta_{x_1}^\eps$
\begin{equation}
\label{eq_geom}
 |M_{x,z}(r'_1,r'_2)| \geq \alpha' (r'_2-r'_1).
\end{equation}
Now, for arbitrary $A_1\subset \Delta_{x_1}^\eps$, $(r'_1,r'_2)\subset ({\hat r}_1,{\hat r}_2)$,
and $x\in \Delta_{x_1}^\eps$, using~(\ref{eq_geom}) we can write
\begin{eqnarray*}
\lefteqn{\PP[\xi^x_2\in A_1,\tau^x_2\in (r'_1,r'_2)]}\\
 &\geq& \int_{\Delta_{y_1}^{\eps'}}dy
                       \int_{A_1}K(x,y)K(y,z)\1{R^x(y,z)\in (r'_1,r'_2)}\, dz \\
  &\geq& h_1^2\alpha' |A_1| (r'_2-r'_1),
\end{eqnarray*}
so $(\xi^x_2,\tau^x_2)$
is $\alpha_1$-continuous on $\Delta_{x_1}^\eps\times ({\hat r}_1,{\hat r}_2)$ with
$\alpha_1 = h_2^2\alpha'$. Lemma~\ref{small_dens} is proved.
\qed

\medskip
We continue the proof of Theorem~\ref{t_conv_KSB}. Fix any $b>1$ such that
\[
b^{-1}{\hat r}_2-b{\hat r}_1 > \max\{(b-1){\hat r}_1,(1-b^{-1}){\hat r}_2\}.
\]
 The following result is a logical continuation
of Lemma~\ref{small_dens}.
\begin{lmm}
\label{l_ind}
For any $k\geq 1$ there exists $\alpha_k>0$ such that for
all $x\in\Delta_{x_1}^\eps$, $(\xi^x_{2k},\tau^x_{2k})$
is $\alpha_k$-continuous on $\Delta_{x_1}^\eps\times (kb{\hat r}_1,kb^{-1}{\hat r}_2)$.
\end{lmm}

\noindent
{\it Proof.}
We proceed by induction. 
The case $k=1$ follows from Lemma~\ref{small_dens}. 
Let $k\geq 2$ and suppose that $(\xi^x_{2(k-1)},\tau^x_{2(k-1)})$
is $\alpha_{k-1}$-continuous on $\Delta_{x_1}^\eps\times ((k-1)b{\hat r}_1,(k-1)b{\hat r}_2)$.
Using Lemma~\ref{small_dens} and the induction assumption, write
for any $A_1\subset \Delta_{x_1}^\eps$, $(r'_1,r'_2)\subset (kb{\hat r}_1,kb^{-1}{\hat r}_2)$
\begin{eqnarray}
\lefteqn{\PP[\xi^x_{2k}\in A_1,\tau^x_{2k}\in (r'_1,r'_2)]}\nonumber\\
  &\geq & \int_{\Delta_{x_1}^\eps}\int_{(k-1)b{\hat r}_1}^{(k-1)b^{-1}{\hat r}_2}
     \PP[\xi^x_{2(k-1)}\in dz,\tau^x_{2(k-1)}\in dt] 
       \PP[\xi^z_2\in A_1,\tau^z_2+t\in (r'_1,r'_2)]\nonumber\\
  &\geq &  \int_{\Delta_{x_1}^\eps}
   \int_{(k-1)b{\hat r}_1}^{(k-1)b^{-1}{\hat r}_2}\alpha_{k-1}\,dz\, dt
    \int_{A_1}\int_{(r'_1-t,r'_2-t)\cap ({\hat r}_1,{\hat r}_2)} \alpha_1 dz'\, dt'\nonumber\\
  &\geq & {\hat C}\alpha_{k-1}\alpha_1|A_1| (r'_2-r'_1), \label{conta_chata}
\end{eqnarray}
where
\[
 {\hat C} = \frac{\min\{(b-1){\hat r}_1,(1-b^{-1}){\hat r}_2\}}
                  {\min\{{\hat r}_2-{\hat r}_1,(k-1)(b^{-1}{\hat r}_2)-b{\hat r}_1\}}.
\]
One can obtain the last inequality in~(\ref{conta_chata}) (together with the value of~${\hat C}$ given
 above) using the following fact:
the convolution of functions $C_1\1{x\in[0,a_1]}$ and $C_2\1{x\in[0,a_2]}$ 
(where $a_1\leq a_2$)
is a function supported on~$[0,a_1+a_2]$, which grows linearly from~$0$ to~$C_1C_2$ on the
interval~$[0,a_1]$, then is flat on the interval~$[a_1,a_2]$, 
then decreases linearly from~$C_1C_2$
to~$0$ on the interval~$[a_2,a_1+a_2]$. So, if one considers any interval $[c,c']$ with
$c>0$ and $c'<a_1+a_2$, then the convolution is uniformly positive 
(and one can easily calculate its minimum) on that interval.
Note that in our case $c=(b-1){\hat r}_1$, $a_1+a_2-c'=(1-b^{-1}){\hat r}_2$.

Lemma~\ref{l_ind} now follows from~(\ref{conta_chata}).
\qed

\medskip

Now, let us choose~$k_0$ in such a way that 
$k_0(b^{-1}{\hat r}_2 - b{\hat r}_1)>(n_0+1)\diam(\DD)$, where~$n_0$ is
 from~(\ref{Doeblin}). Define
\begin{eqnarray*}
 r_3 &=& k_0b{\hat r}_1 + n_0\diam(\DD),\\
 r_4 &=& k_0b^{-1}{\hat r}_2.
\end{eqnarray*}
It is important to observe that $r_4-r_3>\diam(\DD)$.

\begin{lmm}
\label{unif_alpha_cont}
There exists ${\hat\alpha}>0$ such that
$(\xi^y_{n_0+2k_0},\tau^y_{n_0+2k_0})$ is ${\hat\alpha}$-continuous
on $\Delta_{x_1}^\eps\times(r_3,r_4)$ for any $y\in\RR$. 
\end{lmm}

\noindent
{\it Proof.} Consider any $A_1\subset\Delta_{x_1}^\eps$ and
$(r',r'')\subset (r_3,r_4)$. Let $\Phi^y_x(t)$
be the distribution of~$\tau^y_{n_0}$ conditioned on $\{\xi^y_{n_0}=x\}$.
Using~(\ref{Doeblin}) and Lemma~\ref{l_ind} and the fact that
the size of jumps of KRW cannot be larger than~$\diam(\DD)$, write
\begin{eqnarray*}
\lefteqn{\PP[\xi^y_{n_0+2k_0}\in A_1, \tau^y_{n_0+2k_0}\in (r',r'')]}\\
 &\geq& \int_{\Delta_{x_1}^\eps} K^{n_0}(y,x)\,dx\int_0^{n_0\diam(\DD)}
         \PP[\xi^x_{2k_0}\in A_1, \tau^x_{2k_0}\in (r'-t,r''-t)]\, d\Phi^y_x(t)\\ 
 &\geq & {\hat \eps}|\Delta_{x_1}^\eps|\alpha_{k_0}|A_1|(r''-r'),
\end{eqnarray*}
so Lemma~\ref{unif_alpha_cont} holds with 
  ${\hat\alpha}={\hat \eps}|\Delta_{x_1}^\eps|\alpha_{k_0}$.
\qed

\medskip

Now, we are ready to construct a coupling of two versions of the stochastic billiard process
with different starting conditions.
Let $(X_t,V_t)$ be the stochastic billiard with the initial 
condition $(X_0,V_0)=(x_0,v_0)\in\DD\times\Sph^{d-1}$,
and $(X'_t,V'_t)$ has the initial condition $(X'_0,V'_0)=(x'_0,v'_0)$.
Let
\[
 \ttau = \inf\{t\geq 0 : x_0+tv_0\in\fr\}, \qquad w=x_0+\ttau v_0
\]
be the time when $(X_t,V_t)$ first hits the boundary, and the hitting location, and let
\[
 \ttau' = \inf\{t\geq 0 : x'_0+tv'_0\in\fr\}, \qquad w'=x'_0+\ttau' v'_0
\]
be the corresponding quantities for the process $(X'_t,V'_t)$.

 By Lemma~\ref{unif_alpha_cont}, the pairs of 
random variables $(\xi^w_{n_0+2k_0},\ttau+\tau^w_{n_0+2k_0})$ 
and $(\xi^{w'}_{n_0+2k_0},\ttau'+\tau^{w'}_{n_0+2k_0})$ are both ${\hat\alpha}$-continuous 
on $\Delta_{x_1}^\eps\times(r_3,r_4)$. Abbreviate ${\hat h}:=r_4-r_3-\diam(\DD)>0$.
Since $|\ttau-\ttau'|\leq \diam(\DD)$, we obtain that
\[
 |(r_3+\ttau,r_4+\ttau)\cap (r_3+\ttau',r_4+\ttau')| \geq {\hat h}.
\]
So, there exists a coupling of $(\xi^w_{n_0+2k_0},\ttau+\tau^w_{n_0+2k_0})$ and
$(\xi^{w'}_{n_0+2k_0},\ttau'+\tau^{w'}_{n_0+2k_0})$ such that
\begin{equation}
\label{coupl_1st_step}
 \PP[E_1] \geq |\Delta_{x_1}^\eps|{\hat h},
\end{equation}
where
\[
 E_1 = \{\xi^w_{n_0+2k_0} = \xi^{w'}_{n_0+2k_0} \text{ and } 
                    \ttau+\tau^w_{n_0+2k_0} = \ttau'+\tau^{w'}_{n_0+2k_0}\}
\]
(one can use e.g.\ the \emph{maximal coupling}, cf.\ Section~4 of Chapter~1 of~\cite{T}).
On the event~$E_1$, define ${\hat T} = \ttau + \tau^w_{n_0+2k_0}$.

Suppose that the event~$E_1$ does not occur. Consider two cases: 
$\ttau+\tau^w_{n_0+2k_0} \leq \ttau'+\tau^{w'}_{n_0+2k_0}$ and
$\ttau+\tau^w_{n_0+2k_0} > \ttau'+\tau^{w'}_{n_0+2k_0}$. In the first case,
define
\[
 m_1 = \min\{n>n_0+2k_0 : \ttau+\tau^w_n > \ttau'+\tau^{w'}_{n_0+2k_0}\},
\]
and $m'_1=n_0+2k_0$. In the second case, define
\[
 m'_1 = \min\{n>n_0+2k_0 : \ttau'+\tau^{w'}_n > \ttau+\tau^w_{n_0+2k_0}\},
\]
$m_1=n_0+2k_0$ (i.e., we identify the process with the smaller local time, and
let it evolve a few more steps so that its local time becomes close to the other
process' local time). Clearly, in both cases we have
\[
 |(\ttau+\tau^w_{m_1}) - (\ttau'+\tau^{w'}_{m'_1})| \leq \diam(\DD).
\]
Analogously, we obtain that one can couple $(\xi^w_{m_1+n_0+2k_0},\ttau+\tau^w_{m_1+n_0+2k_0})$ and
$(\xi^{w'}_{m'_1+n_0+2k_0},\ttau'+\tau^{w'}_{m'_1+n_0+2k_0})$ in such a way that
\begin{equation}
\label{coupl_2nd_step}
 \PP[E_2\mid E_1^c] \geq |\Delta_{x_1}^\eps|{\hat h},
\end{equation}
where
\[
 E_2 = \{\xi^w_{m_1+n_0+2k_0} = \xi^{w'}_{m'_1+n_0+2k_0} \text{ and } 
                    \ttau+\tau^w_{m_1+n_0+2k_0} = \ttau'+\tau^{w'}_{m'_1+n_0+2k_0}\}.
\]
On the event~$E_1^cE_2$, define ${\hat T} = \ttau + \tau^w_{m_1+n_0+2k_0}$.

Proceeding in this way, we construct two sequences of stopping times $(m_k,m'_k)_{k\geq 1}$
and a sequence of events
\[
 E_k = \{\xi^w_{m_k+n_0+2k_0} = \xi^{w'}_{m'_k+n_0+2k_0} \text{ and } 
                    \ttau+\tau^w_{m_k+n_0+2k_0} = \ttau'+\tau^{w'}_{m'_k+n_0+2k_0}\}
\]
with the property (as in~(\ref{coupl_1st_step}) and~(\ref{coupl_2nd_step}))
\begin{equation}
\label{coupl_kth_step}
 \PP[E_k\mid E_1^c\ldots E_{k-1}^c] \geq |\Delta_{x_1}^\eps|{\hat h};
\end{equation}
on the event $E_1^c\ldots E_{k-1}^cE_k$, we define ${\hat T} = \ttau + \tau^w_{m_k+n_0+2k_0}$.
Since,
\[
 \tau^w_{m_k} - \tau^w_{m_{k-1}} \geq 2(n_0+2k_0)\diam(\DD),
\]
and, by~(\ref{coupl_kth_step})
\[
 \PP[\text{none of the events }E_1,\ldots,E_k\text{ occur}] \leq (1-|\Delta_{x_1}^\eps|{\hat h})^k,
\]
we obtain that, for some $\beta'_1,\beta'_2$
\begin{equation}
\label{coupl_time}
 \PP[{\hat T} > t] \leq \frac{\beta'_1}{2} e^{-\beta'_2 t}.
\end{equation}
Up to time~${\hat T}$, the processes~$X$ and~$X'$ can be explicitly constructed
as in Section~\ref{s_def_model}. Now, by definition we have that at time~${\hat T}$
the local times and the positions are the same for~$X$ and~$X'$, 
and so we have $X_{{\hat T}}=X'_{{\hat T}}$. Since also $X_{{\hat T}}\in\fr$,
we can then construct the two realizations of KSB in such a way that 
$X_{{\hat T}+s}=X'_{{\hat T}+s}$ for all $s\geq 0$
(just use the same sequence of $\eta$s from that moment on). This means that~${\hat T}$
is a coupling time, and, by Theorem~5.1 of Chapter~4 from~\cite{T}, 
we obtain~(\ref{eq_conv_KSB}) from~(\ref{coupl_time}). This 
completes the proof of Theorem~\ref{t_conv_KSB}.
\qed

\medskip
\noindent
{\it Proof of Theorem~\ref{t_cl_KSB}.}
The family of functions $F_t: \bar \DD \times \Sph^{d-1} \to \R$,
$
F_t(x,v)= \int_0^t \E_{x,v} f(X_s,V_s)\, ds 
$
is a Cauchy sequence with respect to the uniform norm $\|\cdot\|_\8$,
in view of (\ref{eq:centered}) and 
(\ref{eq_conv_KSB}). Hence it converges uniformly to a limit~$F$, where
\begin{equation}
  \label{eq:FF}
F(x,v)= \int_0^\8 \E_{x,v} f(X_s,V_s)\, ds.
\end{equation}
Our interest for the limit is that 
\[
M_t=F(X_t,V_t) - \int_0^t f(X_s,V_s)\, ds 
\]
is a martingale, a fact which follows from (\ref{eq:pbmart}) 
since~$F$ satisfies (\ref{eq:dom}) by definition.
We can compute the bracket~$\langle M \rangle$ of the martingale~$M$
(that is, $\langle M \rangle$ is a predictable process, uniquely defined
by the following: $\langle M \rangle_0=0$, and $M_t^2-\langle M \rangle_t$
is a martingale) 
\[
\langle M \rangle_t = 
\sum_{n \geq 1}  {\bf 1}\{\tau_n \leq t\} 
\E_{x,v}\left( [F(X_{\tau_n},V_{\tau_n})-F(X_{\tau_n},V_{\tau_n^-})]^2
\mid \FF_{\tau_n^-}  \right)
\]
Now, with $N(t)\geq 1$ defined by $\tau_{N(t)-1} < t \leq 
\tau_{N(t)}$, we have
\begin{eqnarray*}
\frac{ \langle M \rangle_t}{t} 
&=&
\frac{ \tau_{N(t)}}{t} \frac{N(t)}{ \tau_{N(t)}} \times
 \frac{1}{N(t)}\sum_{n= 1}^{N(t)} 
\E_{x,v}\left( [F(X_{\tau_n},V_{\tau_n})-F(X_{\tau_n},V_{\tau_n^-})]^2
\mid \FF_{\tau_n^-}  \right)\\
& \to& \sigma_f^2 \quad \text{as $t \to \8$}
\end{eqnarray*}
 from Theorem~\ref{t_conv_KRW}, where 
\begin{equation}
\label{eq:varf}
 \sigma_f^2 = \frac{
\E_{\hat \mu} \left( [F(\xi_1,\ell_{\xi_1,\xi_2})-
F(\xi_1,\ell_{\xi_0,\xi_1})]^2  
 \right)}{
\E_{\hat \mu} \|\xi_1-\xi_0\|}
\in [0,\8) .
\end{equation}
By the central limit theorem for martingales 
\cite[theorem 35.12]{Bill},
it follows that 
\[
 t^{-1/2} M_t \to  \text{\rm Normal}(0,\sigma_f^2)
\]  
in law. Since~$F$ is bounded, we see that 
\[
t^{-1/2} \int_0^t f(X_s,V_s) ds = - t^{-1/2} M_t + \OO(t^{-1/2}),
\]
which implies the desired convergence and concludes the proof of Theorem~\ref{t_cl_KSB}.

Let us now prove that~(\ref{func_var0}) implies that $\sigma^2_f=0$.
For KSB, this function~$f$ is centered with respect to the uniform measure~$\chi$ 
which is invariant. 
In this case, we can compute the function $F$ from~(\ref{eq:FF}):
\begin{eqnarray*}
F(x,v) &=&
- \int_0^\8 \E_{x,v} V_s \cdot \nabla G(X_s) ds\\
&=& -
\lim_{t \to \8} \E_{x,v} \int_0^t  V_s \cdot \nabla G(X_s) ds\\
&=&
G(x)- \lim_{t \to \8} \E_{x,v}  G(X_t) \\
&=& G(x)  - \int_{\DD \times \Sph^{d-1}} G(x) \chi(dx,dv) 
\end{eqnarray*}
which does not depend on $v$. Hence  $F(x,\cdot)$ is constant on 
$\Sph_{\n(x)}$, and the definition (\ref{eq:varf}) shows that
  $\sigma_f^2=0$. This can be understood in a different manner
by computing 
\begin{eqnarray*}
M_t &=&  G(X_t) + \int_0^t   V_s \cdot \nabla G(X_s) ds\\
 &=&  G(x) - \int_{\DD \times \Sph^{d-1}} G(x) \chi(dx,dv)  
\end{eqnarray*}
Therefore,
\[
 \int_0^t f(X_s,V_s) ds = G(X_t)-G(x)
\]
is bounded, which shows again that $\sigma_f^2=0$.

It is not clear to us if one can find a condition on~$f$ which is necessary and sufficient
for $\sigma_f^2=0$. Notice however that if for a.e.\ $x \in \fr$, 
$F(x,\cdot)$ is constant on $\Sph_{x}$, then $\sigma_f^2=0$.
\qed

\medskip
\noindent
{\it Proof of Theorem~\ref{th_gen}.} 
We identify the invariant measure $\chi$ appearing in~(\ref{eq_conv_KSB}).
By theorems \ref{t_conv_KSB} and  \ref{t_conv_KRW}, we have for bounded 
continuous functions $f: \bar \DD \times \Sph^{d-1} \to \R$ and starting from 
any initial condition,
\begin{eqnarray}
\int_{ \bar \DD \times \Sph^{d-1}} f(x,v) \,dx\, dv 
&=&
\lim_{t \to \8} t^{-1} \int_0^t \E f(X_s, V_s) \,ds  \nonumber \\
&=&
\lim_{n \to \8} \E \tau_n^{-1} \sum_{k=1}^n \int_{\tau_{k-1}}^{\tau_{k}}  
f(X_s, V_s) \,ds  \nonumber \\
&=&
 \big({\E_{\hat \mu} \|\xi_1-\xi_0\|}\big) ^{-1} 
\E_{\hat \mu} 
\int_{0}^{\tau_1}  f(\xi_0+ s\ell_{\xi_0,\xi_1})\,ds.\phantom{*} \label{eq:chi1}
\end{eqnarray}
The numerator in (\ref{eq:chi1}) is equal to 
\begin{eqnarray*}
\lefteqn{\int_{\partial \DD} \hat \mu(dz) \int_{\Sph_{\n(z)}} dv \;
\bga( U_z^{-1} v) \int_0^\8 ds f(z+sv,v) \1{z \leftrightarrow z+sv}}\\
&=&\int_{\DD}dx \int_{\partial \DD} dz {\|z-x\|^{-d+1}} \frac{d \hat \mu}{dz}(z)
\bga( U_z^{-1} \ell_{z,x}) f(x,\ell_{z,x}) \1{z \leftrightarrow x} 
\end{eqnarray*}
with the change from polar coordinates $(s,v)$ around~$z$ to Cartesian
coordinates~$x$. Then, by performing the change of variables from
$z \in \fr$ to $v=\ell_{z,x} \in  \Sph^{d-1}$, the Hausdorff measure integrator
$dz {\|z-x\|^{-d+1}}$ becomes  the Haar measure integrator
$dv/\cos \phi_z(v)$ with $z={\mathsf h}_{x}(-v)$, and the
whole integral becomes
\[
\int_{\DD}dx \int_{ \Sph^{d-1} } dv f(x,v) 
\frac{d \hat \mu}{dz}({{\mathsf h}_{x}(-v)})
\frac{\bga( 
U^{-1}_{{\mathsf h}_{x}(-v)} v) }{\cos \phi_{{\mathsf h}_{x}(-v)}(v)}.
\]
This is indeed the formula claimed in  Theorem~\ref{th_gen}.
\qed

\medskip
\noindent
{\it Proof of Theorem~\ref{th_gen2}.} 
Since the process  $(\tX_t,\tV_t)_{t \in \R}$ is stationary and Markov,
it is enough to show that 
$(\tX_0,\tV_0,\tX_t,\tV_t) \eqlaw (\tX_t,-\tV_t,\tX_0,\tV_0)$ 
for all $t \geq 0$, i.e., that
\begin{equation} \label{eq:quasirevf}
\E f(\tX_0,\tV_0,\tX_t,\tV_t)= \E f(\tX_t,-\tV_t,\tX_0,-\tV_0)
\end{equation}
for all smooth test functions  $f$  on
$ \DD \times \Sph^{d-1} \times \DD \times \Sph^{d-1}$.
Since
\[
\E f(\tX_0,\tV_0,\tX_t,\tV_t)=
\sum_{n \geq 0} \E \left[f(\tX_0,\tV_0,\tX_t,\tV_t) \1{\tau_n \leq t < 
\tau_{n+1}} \right]
\]
this equality (\ref{eq:quasirevf}) will follow from the relation
\begin{eqnarray} 
 \lefteqn{\E \Big[f(\tX_0,\tV_0,\tX_t,\tV_t) \1{\tau_n \leq t < 
\tau_{n+1}} \Big]} \nonumber\\
&=&
 \E \Big[f(\tX_t,-\tV_t,\tX_0,-\tV_0) \1{\tau_n \leq t < 
\tau_{n+1}} \Big]\label{eq:quasirevfn}
\end{eqnarray}
for $n=0,1,\ldots$. The case $n=0$ 
being clear, we now consider the case 
$n \geq 1$.
The left-hand side of (\ref{eq:quasirevfn}) is equal to
\[
 \int_{A_n}
 dx\, dv\, dy_2 \ldots dy_{n+1} 
 \tK(y_1,y_2)  \ldots \tK(y_n,y_{n+1}) f(x,v, y_n+R \ell_{y_n,y_{n+1}}, \ell_{y_n,y_{n+1}})
\]
where $y_1= {\mathsf h}_x(v)$, where $R$ is the function of
$x,y_1,\ldots,y_n$ given by
\[
R=
t-s_{1,n}-\|y_1-x\|\;,
\quad s_{k,n}=\|y_{k+1}-y_k\|+ \ldots +
\|y_{n}-y_{n-1}\|\;,
\]
and $A_n \subset \DD \times \Sph^{d-1}  \times (\fr)^{n}$ is defined
by $\|y_1-x\|+s_{1,n} \leq t < \|y_1-x\|+s_{1,n+1}$. 
Changing the Cartesian variable $x$ into the polar coordinates 
$s,w$ around $y_1$,
$s \in [0, \|y_1-y_0\|], w \in \Sph^{d-1}$ ($y_0$ being defined as
$y_0= {\mathsf h}_x(-v)$), we write
\begin{eqnarray*}
dv\,dx
&=&
dv  \|y_1-y_0\|^{d-1} \1{s\leq \|y_1-y_0\|}  dw\, ds\\
&=&
dw  \frac{\|y_1-y_0\|^{d-1}}{\cos \phi_{y_1}(w)} 
dv  \frac{\|y_1-y_0\|^{d-1}}{\cos \phi_{y_0}(v)}
\tK(y_0,y_1) \1{s\leq\|y_1-y_0\|}  ds\\
&=&
dy_0 \, dy_1 \tK(y_0,y_1) \1{s\leq\|y_1-y_0\|} ds
\end{eqnarray*}
in terms of the Hausdorff measure on~$\fr$. Finally, the left-hand side 
of (\ref{eq:quasirevfn}) writes
\begin{eqnarray*}
\lefteqn{\int_{R \leq \|y_{n+1}-y_n\|, s \leq \|y_1-y_0\|}
dy_0  \ldots dy_{n+1} \, ds
\tK(y_0,y_1)  \ldots \tK(y_n,y_{n+1})}\\
&&\qquad \qquad 
\times f\left( y_1-s  \ell_{y_0,y_1} ,  \ell_{y_0,y_1},
 y_n+R  \ell_{y_n,y_{n+1}}, \ell_{y_n,y_{n+1}}\right)
\end{eqnarray*}
with~$R$ defined by $R+s+s_{1,n}=t$. Defining now  $S$ by $S+r+s_{1,n}=t$,
and changing variables to $r=R, z_0=y_{n+1}, \ldots z_{n+1}=y_0$, the 
integral becomes
\begin{eqnarray*}
\lefteqn{\int_{S \leq \|z_{n+1}-z_n\|, r \leq \|z_1-z_0\|}
dy_0  \ldots dy_{n+1}\, ds
\tK(z_1,z_0)  \ldots \tK(z_{n+1},z_n)}\\
&&\qquad \qquad 
\times f\left( z_n+S \ell_{z_n,z_{n+1}},-\ell_{z_n,z_{n+1}},  \ell_{y_0,y_1},
 z_1-r  \ell_{y_0,y_1}  ,- \ell_{y_0,y_1} \right)
\end{eqnarray*}
which is the right-hand side of (\ref{eq:quasirevfn}) since
$\tK$ is a symmetric function. This ends the proof. 
\qed

\section{Proofs of geometric properties of the random chord}
\label{s_proofs_RC}

\medskip
\noindent
{\it Proof of Theorem~\ref{th_chord}.}
To keep the paper self-contained, we provide the proof here, even though
there are a lot of similarities with the proof of corresponding results
in~\cite{B1,E,LR}.

Let us start the process from the uniform distribution on~$\fr$ (recall that
it is the stationary distribution for KRW). Using~(\ref{eq_conv_KSB})
and the ergodic theorem one can write that
\begin{eqnarray*}
\frac{|A|}{|\DD|} &=& \lim_{t\to\infty} t^{-1} \int_0^t\1{\tX_s\in A}\, ds\\
&=& \lim_{n\to\infty} \tau_n^{-1}\int_0^{\tau_n} \1{\tX_s\in A}\, ds \\
&=& \lim_{n\to\infty} \Big(n^{-1}\sum_{i=1}^n(\tau_i-\tau_{i-1})\Big)^{-1} 
   n^{-1}\sum_{i=1}^n \int_{\tau_{i-1}}^{\tau_i} \1{\tX_s\in A}\, ds\\
&=& \bm^{-1}\m(A),
\end{eqnarray*}
so we already obtain that $\m(A)$ is proportional to~$|A|$ (with a factor that
does not depend on~$A$):
\begin{equation}
\label{proportional_to_|A|}
\m(A) = \frac{\bm}{|\DD|}|A|.
\end{equation}

To calculate the value of the factor $\frac{\bm}{|\DD|}$, take any~$x\in\DD$ and a small~$\delta>0$
(so that $\BB(x,\delta)\subset\DD$). For any $u\in\Sph^{d-1}$ denote by
\[
 G_u = \{y\in\fr : \text{ there is }s>0 \text{ such that }y+su\in\BB(x,\delta)
         \text{ and }y\leftrightarrow y+su\}
\]
the projection of $\BB(x,\delta)$ in the $(-u)$-direction. For $y\in\R^d$ let
\[
 H_{y,u} = |\{s:y+su\in\BB(x,\delta)\}|
\]
be the length of the intersection of the ray drawn from~$y$ in the $u$-direction
with $\BB(x,\delta)$. Then
\begin{eqnarray}
\m(\BB(x,\delta)) &=& \frac{1}{|\fr|}\int_{\fr}dy
        \int_{\Sph_{\n(y)}}\gamma_d\cos\phi_y(u)H_{y,u}\,du\nonumber\\
 &=& \frac{\gamma_d}{|\fr|}\int_{\Sph^{d-1}}du
                \int_{G_u}\cos\phi_y(u)H_{y,u}\,dy.
 \label{change_order_int}
\end{eqnarray}
For almost all~$u\in\Sph^{d-1}$ there exists an 
unique $y_u\in\fr$ that ``sees~$x$ in
direction~$u$'', i.e., $\{y_u\leftrightarrow x\}$ 
and there exists~$t_u>0$ such that $y_u+t_u u = x$. Define by
\[
 {\hat G}_u = \{z\in\R^d : (z-y_u)\cdot\ell_{y_u,x}=0, 
        \text{ there is }s>0 \text{ such that }y+su\in\BB(x,\delta)\}
\]
the projection of $\BB(x,\delta)$ onto the hyperplane that passes
through~$y_u$ and is perpendicular to $\ell_{y_u,x}$. Then,
provided that $\cos\phi_{y_u}(u)>0$,
\begin{eqnarray*}
\int_{G_u}\cos\phi_y(u)H_{y,u}\,dy &=& \int_{G_u}(\cos\phi_{y_u}(u)+O(\delta))H_{y,u}\,dy\\
 &=&\int_{{\hat G}_u} (1+O(\delta)) H_{z,u}\, dz \\
 &=& |\BB(x,\delta)|(1+O(\delta)),
\end{eqnarray*}
so (\ref{proportional_to_|A|}) and (\ref{change_order_int}) imply that
\[
 \frac{\bm}{|\DD|} = \frac{\gamma_d|\Sph^{d-1}|}{|\fr|},
\]
and we finally obtain~(\ref{mean_chord}) from~(\ref{proportional_to_|A|}). \qed

\medskip
\noindent
{\it Proof of Theorem~\ref{t_restr_convex}.}
First, we observe that the boundary of any convex domain is
almost everywhere continuously differentiable. Indeed, since the boundary of a convex domain 
is locally a graph of a convex function, this follows from Theorems~D and~E of Section~44 
of~\cite{RV}.

\begin{figure}
\centering
\includegraphics{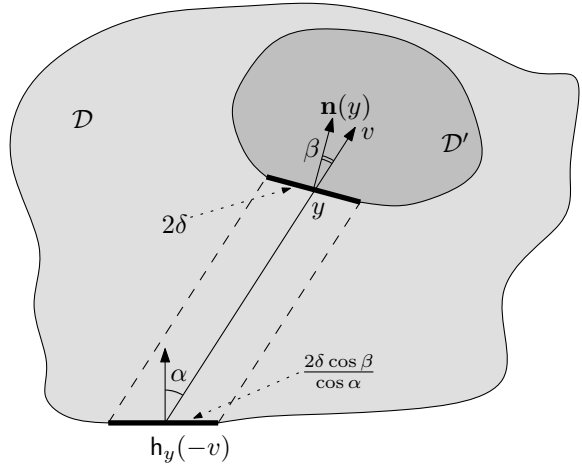}
\caption{On the proof of Theorem~\ref{t_restr_convex}}
\label{f_restr_convex}
\end{figure}

Consider a point $y\in\fr'$ where $\fr'$ has a locally $\cCC^1$
parametrization.
As before, denote by
\[
 G'_u = \{z\in\fr : \text{ there is }s>0 \text{ such that }z+su\in\fr'\cap\BB(y,\delta),
         z\leftrightarrow z+su\}
\]
the projection of $\fr'\cap\BB(y,\delta)$ onto~$\fr$ in the $(-u)$-direction. 
Let $I=\{[\Xi_1,\Xi_2] \cap \DD' \neq \emptyset \}$ be the event that the 
random $\DD$-chord intersects $\DD'$.
For $B\subset \Sph_{\n(y)}$, write (cf.\ Figure~\ref{f_restr_convex}, and
note that $\ell_{\Xi_2,\Xi_1}=\ell_{\Xi'_2,\Xi'_1}$)
\begin{eqnarray*}
\lefteqn{\PP[\Xi'_1\in\BB(y,\delta), \ell_{\Xi'_2,\Xi'_1}\in B]}&&\\ 
&=& \PP[\Xi_1\in G'_{\ell_{\Xi_2,\Xi_1}}, \ell_{\Xi_2,\Xi_1}\in B
                \mid I ] \\
     &=& \PP[\Xi_1\in G'_{\ell_{\Xi_2,\Xi_1}}, \ell_{\Xi_2,\Xi_1}\in B] \times 
                   \PP(I)^{-1} \\
     &=& \frac{\gamma_d}{|\fr|\;\PP(I)}\int_B dv \int_{G'_v}\cos\phi_x(v)\, dx\\
&=& \frac{\gamma_d}{|\fr|\;\PP(I)} (|\fr'\cap \BB(y,\delta)|+
o(\delta^{d-1})) \int_B \cos\phi_y(v) \, dv.  
\end{eqnarray*}
This shows in fact that the couple $( \Xi'_1, \ell_{\Xi'_2,\Xi'_1})$ 
has a density
with respect to $dy\, dv$ which is proportional to $\cos\phi_y(v) 
\1{\Sph_{\n(y)}}$, i.e.,
\begin{equation}  
\label{eq:rocky}
\PP[\Xi'_1\in\BB(y,\delta), \ell_{\Xi'_2,\Xi'_1}\in B] =
\frac{\gamma_d}{|\fr|\;\PP(I)} \int_{(\fr' \cap \BB(y,\delta))\times
B} \cos\phi_x(v) \, dx\,dv,
\end{equation}
which proves Theorem~\ref{t_restr_convex}.
\qed

\medskip
\noindent
{\bf Remark.} We can also identify the normalization in the last formula,
\begin{equation}
  \label{eq:rocky2}
\PP(I)= \frac{|\fr'|}{|\fr|}.
\end{equation}

\medskip
\noindent
{\it Proof of Theorem~\ref{t_nonconvex}.}
Let  $y\in\fr'$ where $\fr'$ has a locally $\cCC^1$
parametrization, $\delta>0$ 
and a Borel set $B\subset \Sph_{\n(y)}$. Similarly to~(\ref{eq:rocky}), we get
\begin{multline}
\PP\Big[\exists k \leq \iota: \Xi'_{1,k}\in\BB(y,\delta), 
\ell_{\Xi'_{2,k},\Xi'_{1,k}}\in B\Big] \\
=
 \frac{\gamma_d}{|\fr|} 
\int_{(\fr' \cap \BB(y,\delta))\times
B} \cos\phi_x(v) \, dx\,dv.
\end{multline}
Note that left-hand side can be written as
\[
 \E \Big[\sum_{k =1}^{ \iota } \1{\Xi'_{1,k}\in\BB(y,\delta), 
\ell_{\Xi'_{2,k},\Xi'_{1,k}}\in B}\Big].
\]
Then, for  $A_1, A_2\in \partial\DD'$ such that $A_1\cap A_2=\emptyset$ and 
$x \leftrightarrow y$ in $\DD'$ for all $x \in A_1, y \in A_2$, we get
\[
\E \Big[ \sum_{k =1}^{ \iota} \1{\Xi_{1,k}'\in A,\,\Xi_{2,k}'\in B}
 \Big] = \frac{1}{|\fr|}
\int_{A\times B}\tilde K(x,y)\,dx\, dy.
\]
This is enough to ensure that, for  
 $C \subset \partial\DD'\times \partial\DD'$, we have
\begin{equation*} 
\E \Big[ \sum_{k =1}^{ \iota} \1{(\Xi_{1,k}',\Xi_{2,k}')\in C}
\Big] =  \frac{1}{|\fr|}
\int_{C}\tilde K(x,y)\,dx\, dy.
\end{equation*}
which proves (\ref{eq:baila}). Taking
$C=\partial\DD'\times \partial\DD'$ in (\ref{eq:baila}),
 we obtain~(\ref{eq:baila1}). \qed

\medskip
\noindent
{\it Proof of Theorem~\ref{t_intersection}.}
For any $v\in\Sph_e$ and $\delta>0$ define
\[
 S_v(\delta) = \{z\in\R^d : \text{ there exist }t\in[-\delta,\delta], x\in S
                 \text{ such that }z+tU_x v = x\}.
\]
Note that
\begin{equation}
\label{vol_Sdelta}
 |S_v(\delta)| = 2\delta (1+o(\delta)) |S| \cos\phi(e,u).
\end{equation}

Let us consider first the process in the stationary regime, i.e., 
we suppose that $\tX_0,\tV_0$ (and hence $\tX_t,\tV_t$) are
independent and uniform. For $\delta<\inf_{x\in\fr,y\in S}\|x-y\|$
write (observe that the next probability
does not depend on~$t$)
\begin{eqnarray*}
\lefteqn{\PP[\text{there exists $n$ such that }\htau_n\in[t,t+\delta],w_n\in B]}\\
&=& \PP[\tX_t\in S_{\tV_t}(\delta)]\\
&=& \frac{1}{|\Sph^{d-1}|}\int_B \PP[\tX_t\in S_u(\delta)\mid \tV_t=u]\, du\\
&=& \frac{1}{|\Sph^{d-1}|}\int_B \PP[\tX_t\in S_u(\delta)]\, du\\
&=& \delta (1+o(\delta))|S| \frac{2}{|\Sph^{d-1}|}\int_B \cos\phi(e,u)\, du.
\end{eqnarray*}
This implies~(\ref{eq_intersection}) for the stationary case. For the general case,
the formula~(\ref{eq_intersection}) now easily follows from Theorem~\ref{t_conv_KSB}.
\qed

\medskip
\noindent
{\bf Acknowledgements:} We thank Fran\c cois Delarue and Stefano Olla for
stimulating discussions.

\end{document}